\documentclass[onefignum,onetabnum]{siamart190516}

\usepackage{tensor}
\usepackage{lipsum}
\usepackage{amsfonts}
\usepackage{graphicx}
\usepackage{epstopdf}
\usepackage{algorithmic}
\ifpdf
  \DeclareGraphicsExtensions{.eps,.pdf,.png,.jpg}
\else
  \DeclareGraphicsExtensions{.eps}
\fi

\newsiamremark{remark}{Remark}
\newsiamremark{hypothesis}{Hypothesis}
\crefname{hypothesis}{Hypothesis}{Hypotheses}
\newsiamthm{claim}{Claim}

\headers{}{Lilia Ghaffour, and Taous-Meriem Laleg-Kirati}

\title{Reference Tracking \R{AND} Observer Design for Space\R{-}Fractional Partial Differential Equation Modeling Gas Pressure\R{s} in Fractured Media{%Submitted to the editors DATE.
}}

\author{Lilia Ghaffour \thanks{Computer, Electrical and  Mathematical Science and  Engineering division at 
	King Abdullah University of Science and Technology, Thuwal 23955-6900, Saudi Arabia
  (\email{lilia.ghaffour@kaust.edu.sa}, \email{taousmeriem.laleg@kaust.edu.sa} ).}
\and  Taous-Meriem Laleg-Kirati \footnotemark[1] \thanks{ National Institute for Research in Digital Science and Technology (INRIA), France
  .}
}

\usepackage{amsopn}

\usepackage{caption}
\usepackage{subcaption}
\newcommand{\B}[1]{\textcolor{black}{#1}}
\newcommand{\R}[1]{\textcolor{black}{#1}}

\newcommand\norm[1]{\left\lVert#1\right\rVert}
\usepackage{lipsum}
\usepackage{amsfonts}
\usepackage{graphicx}
\usepackage{epstopdf}
\usepackage{algorithmic}
\ifpdf
  \DeclareGraphicsExtensions{.eps,.pdf,.png,.jpg}
\else
  \DeclareGraphicsExtensions{.eps}
\fi

\begin{document}

\maketitle

% REQUIRED
\begin{abstract}
 This paper \R{considers} a class of space fractional partial differential equations \B{(FPDEs)} \R{that} describ\R{e} gas pressure\R{s} in fractured media. First, the well-posedness, uniqueness, and the stability in $\B{L_{\infty}(\mathbb{R})}$of the \R{considered FPDEs} are investigated. Then, the reference tracking \B{problem is studied to track} the pressure gradient at a downstream location of a channel. \R{This requires manipulation of} gas pressure at \R{the} downstream location and \R{the use of pressure} measurements at an upstream location. \B{To achiever this, the backstepping approach is adapted to the space FPDEs.} The key challenge in this adaptation is the non\R{-}applicability of the Lyapunov theory which is \R{typically} used to prove the stability of the \B{target system as, the obtained target system is fractional in space.}  In addition, a \B{backstepping adaptive observer is  designed  to jointly estimate both the system's state and the disturbance}. \B{The stability of the closed loop (reference tracking controller/observer) is \R{also} investigated. Finally, numerical simulations are given to evaluate the efficiency of the proposed method.}
\end{abstract}

% REQUIRED
\hspace{-3mm}\begin{keywords}
 \hspace{-0.8mm}Reference Tracking,\hspace{-0.2mm} Observer Design, \hspace{-0.4mm}Fractional systems, \hspace{-0.4mm}Well-posedness\R{, \hspace{-0.4mm}U}niqueness\R{.}
\end{keywords}

% REQUIRED
\section{INTRODUCTION}
Asymptotic output regulation, \R{which} allow\R{s} a system's output to follow a desired reference, is  a challenging problem \R{that} has \R{received considerable research} attention\R{, particularly with regard to the }presence of disturbances. This problem has been widely studied and solved for different perturbed systems  described \R{either using} \R{o}rdinary \R{d}ifferential \R{e}quations (ODEs) or \R{p}artial \R{d}ifferential \R{e}quations (PDEs), including parabolic PDEs,
which describe a wide variety of time-dependent phenomena, such as heat conduction and particle diffusion \R{(e.g., \cite{meurer2012control}, \cite{deutscher2015backstepping})}. For \R{example}, in \cite{isidori1993topics}, \R{the} output regulation problem was investigated for lumped\R{-}parameter linear systems. An extension of this result was discussed in \cite{byrnes1997output} \R{and} \cite{isidori2012robust} for nonlinear systems. \R{Subsequent} studies \R{have} investigated the robust output regulation problem  for linear distributed parameter systems \R{(e.g., \cite{bymes2000output}, \cite{immonen2007internal}, \cite{hamalainen2010robust})}. More recently, \R{this approach has been} extended to different PDE systems  \cite{paunonen2013output}, \cite{natarajan2014state}. In particular, in  \cite{deutscher2015backstepping}, the author considered \R{the} output regulation problem for boundary\R{-}controlled parabolic PDEs in the presence of disturbances under the assumption that these disturbances and the desired reference can be modeled \R{using a} known finite dimensional signal. 

\B{On the other hand, 
fractional PDEs} \B{(FPDEs)} have been used for accurately  modeling  and analyzing numerous  phenomena in different scientific and engineering fields \R{(e.g., \cite{meerschaert2004finite}, \cite{xiong2012inverse}) owing} to their suitability to describe \R{the dynamics of} several problems. Both time and space FPDEs (fractional derivative\R{s} in time and space\R{,} respectively) \R{can} capture multi scale features of complex physical phenomena \R{due} to the non locality of the fractional operators \cite{basu2002quadratic},  \cite{benson2000application}, \cite{benson1998fractional}, \cite{el2002continuation} and \cite{mainardi1997fractional}. 
\B{Furthermore,} time FPDEs have been more \R{widely investigated} in the literature than space FPDEs\R{, predominantly in the field of control.}  For example,   \cite{liang2004boundary} propose\R{d} to solve boundary feedback stabilization \R{and disturbance
rejection problems for time fractional PDE. whereas}  \cite{zhou2018boundary} solves the  output feedback \R{stabilization problem for an unstable time fractional PDE}. \R{ \cite{li2012observer} and \cite{song2013dynamic} addresses the observer-based robust stabilisation problem and dynamic output feedback control for non-linear fractional uncertain systems and fractional order systems} respectively. \cite{fujishiro2014approximate} \R{considers approximate  controllability problem for fractional  diffusion  equations.}

In \cite{schumer2009fractional}, the authors concluded that particle transport behavior may be parsimoniously described \R{using} a fractional advection dispersion equation (FADE) by phenomenological discussion \R{of the} arbitrary average velocity and non\R{-}zero dispersion coefficient. 
%In \cite{zheng2010spectral}, the authors considered a Cauchy problem for time FADE, they presented a spectral regularization method. 
\R{Furthermore}, in  \cite{chen2008adi}, an unconditionally stable second\R{-}order difference method for two\R{-}dimensional FADE \R{was} proposed. 

\B{Moreover, }very few work investigates the tracking problem for FPDEs. \B{\cite{ghaffour_laleg_kirati} investigates the reference tracking problem for time fractional advection dispersion equation in presence of disturbances}.
Less investigations have been done on space fractional PDEs despite their importance in describing a large variety of real life phenomena such as gas production in fractured media \cite{amir2018physics}, solute transport in heterogeneous porous media \cite{schumer2009fractional}, \cite{schumer2001eulerian}, Plumes spread in laboratories in \cite{benson2000application} \cite{benson1998fractional}, and  Transport affected by hydraulic conditions at a distance on the Earth surface in \cite{chen2008adi}. To the best of the authors knowledge, there is no paper dealing with stability,  reference tracking and observer design  problems for space FPDE. \B{Furthermore,} some work exists on studying the well-posedness of space FDEs. For instance, in  \cite{huang2005fundamental}, the authors solved analytically and in terms of Green function the homogeneous time and space FADE. and  in  \cite{aldoghaither2017direct}, the non homogeneous Riesz-Feller space FADE has been solved  analytically.

In this paper, we consider a class of Caputo space FPDE, which describes for instance the gas production in fractured media. such mechanisms is done by drilling wells through gas saturated rock to force 
the gas to flow through the drilled well into the production pipelines. The underground layers burden pressure causes a chaotic explosion of the saturated rocks,  which yields to unequal spatial distribution of the pressure and which is modeled by a spatial fractional derivative. In \cite{amir2018physics}, the authors solves the non-Newtonian non-Darcy fractional-derivatives flux equations using physics-preserving averaging schemes that incorporates both, original and shifted, Grunwald- Letnikov (GL) approximation formulas preserving the physics, by reducing the shifting effects, while maintaining the stability of the system. They derived the system's  equations and discussed the discretization schemes. Then, they illustrated the physics-preserving averaging scheme. Some authors believe that a minimum pressure gradient (called threshold pressure gradient (TPG)) is required before a liquid starts to flow in a porous medium\B{. I}t has been proven as well in \cite{sheikha2009effect} that the pressure gradient has a much greater effect on gas mobility and oil recovery than pressure-decline rate has. That is why, in this paper, we aim at tracking the pressure gradient at a final position to follow a desired trajectory. This tracking process can be seen as a decreasing process of the gradient pressure as well by taking a desired trajectory which is as small as possible.

The main objectives of this paper are to study the  reference tracking problem for boundary controlled space FPDE \B{modelling the gas production in fractured media} in the presence of disturbances. \B{This is done by adapting the backstepping approach to the space FPDEs.} The key challenge here is the non applicability of the Lyapunov theory which is commonly used to prove the stability of the target system \B{as the obtained target system is fractional in space}. \B{The contribution of this paper is to prove the asymptotic stability of the space fractional target system to ensure the reference tracking process. This is done by deriving the analytical solution of the target system in terms of Green functions. Then, taking the advantage of the fractional differentiation order which will allow the green function to converge asymptotically to zero.}
\B{With a similar reasoning another contribution of the paper would be to extract some conditions on the system's parameters that ensure the stability \R{ in $\B{L_{\infty}(\mathbb{R})}$} of the considered model .} Moreover, \B{an additional contribution of this paper is the design of a backstepping adaptive state observer which is required by the output feedback control}. \B{In addition, the stability of the closed loop (reference tracking controller-observer) is proved.} Moreover, a fundamental solution  for the non homogeneous Caputo space FPDE is proposed and its uniqueness is studied. 

This paper is organized as follows, in section 2, a review of some preliminaries on fractional calculus is presented. In section 3, the considered problem is formulated. The fundamental solution of the considered FPDE using Duhamel's principle and the generalized Leibnitz differentiation rule is introduced in section 4. Section 5 aims to derive some sufficient conditions to guarantee the stability \R{in $\B{L_{\infty}(\mathbb{R})}$} of the considered problem. Section 6 studies the reference tracking problem for boundary controlled space FPDE in the presence of disturbances by proving the asymptotic stability \R{in $\B{L_{\infty}(\mathbb{R})}$} of the tracking error. The \B{backstepping adaptive} observer \B{design for} the considered problem is \B{given} in section 7. \B{Section 8 studies the stability of the closed loop (reference tracking controller-observer). Simulation examples are given as well to support the proposed results.} Finally, a general conclusion will summarize the obtained results.

 \section{Fractional Calculus Preliminaries}
This section presents useful  definitions and results on fractional PDE. 
\begin{definition}[\R{Gamma function}]
\cite{podlubny1998fractional, ghaffour2018fractional, ghaffour2016class}
\R{Gamma function is defined by:
\begin{equation*} 
\Gamma(z)=\int_{0}^{+\infty} e^{-t}t^{z-1}dt, \ \quad z \in \mathbb{R^{+}}, 
\end{equation*} }
\end{definition}
\begin{definition} {\cite{podlubny1998fractional,ghaffour2018fractional, ghaffour2016class}}\\
Riemann-Liouville fractional integral of order $\alpha$ with $n-1 < \alpha \leq n, n \in {\rm I\!N}$ is defined by:\\
\vspace{-1cm}
\begin{equation*}
_{a}\mathcal{J}^{\alpha}_{t}f(t)=\frac{1}{\Gamma(\alpha)} \int _{a} ^{t} (t-\tau)^{\alpha-1} f(\tau)d\tau,   
\\
\end{equation*}
where,  $t\in [a,b], f(.) \in C[a,b]$. We denote:
$\mathcal{J}^{\alpha}_{t}f(t):=  \quad _{0}\mathcal{J}^{\alpha}_{t}f(t). $ 
\\
\end{definition}

\begin{definition}{\cite{podlubny1998fractional,ghaffour2018fractional, ghaffour2016class} }\\
Caputo fractional derivative of order $\alpha$ with $n-1 < \alpha \leq n, n \in {\rm I\!N}$ is given by:\\
\begin{equation*} 
^{C}_{a}D^{\alpha}_{t} f(t)=
\left\{\begin{array}{llll}    
\frac{\partial ^{n}f(t)}{\partial t^{n}},  \quad  \alpha = n \in {\rm I\!N}
\\
\frac{1}{\Gamma(n-\alpha)} \int _{a} ^{t} (t-\tau)^{n-\alpha-1} \frac{\partial ^{n}f(\tau)}{\partial \tau ^{n}}d\tau,
\end{array}\right.
\end{equation*}
where,  $t\in [a,b], f(.) \in C^{n}[a,b]$. We denote:
$$ ^{C}D^{\alpha}_{t} f(t):= \quad ^{C}_{0}D^{\alpha}_{t} f(t). $$ 

For more information on the Caputo fractional derivative, we refer the readers to \cite{podlubny1998fractional}.
\end{definition}

\begin{remark}
Caputo time fractional derivative with a negative non integer order is defined in \cite{podlubny1998fractional} by:
\begin{equation*} 
^{C}D^{\alpha}_{t}f(t)=\mathcal{J}^{-\alpha}_{t}f(t),
\end{equation*}

for all $\alpha \in \mathbb{R}^{-}/\mathbb{Z}^{-}$.
\end{remark} 
\begin{theorem}
Let $f(t) \in C^{n}[a,b], t\in [a,b],$ we have:
\begin{equation*}
^{C}D^{\alpha}_{t} \mathcal{J}^{\alpha}_{t} f(t)=f(t),
\end{equation*} 
where $n-1 < \alpha \leq n, n \in {\rm I\!N}$, see \cite{podlubny1998fractional} for the proof.
\end{theorem}

\begin{definition} { \cite{podlubny2007adjoint} }\\
The right-sided Caputo fractional derivative of order $\alpha$ with $n-1 < \alpha \leq n, n \in {\rm I\!N}$ is defined by:
\vspace{-0.5cm}
\begin{equation}\label{sided} 
\tensor*[^{C}]{D}{^{\alpha}_{t,b}}f(t)=\frac{(-1)^{n}}{\Gamma(n-\alpha)} \int _{t} ^{b} (\tau-t)^{n-\alpha-1} \frac{\partial ^{n}f(\tau)}{\partial \tau ^{n}}d\tau, 
\\
\end{equation}
where,  $t\in [a,b], f(.) \in C^{n}[a,b]$.
\end{definition}

\begin{definition} { \cite{podlubny2007adjoint} }\\
\R{The right-sided Caputo time fractional derivative with a negative non integer order $\alpha \in \mathbb{R}^{-}/\mathbb{Z}^{-}$ is defined by:}
\R{ \begin{equation*} 
^{C}D^{\alpha}_{t,b}f(t)=\mathcal{J}^{-\alpha}_{t,b}f(t),
\end{equation*}
where, $\mathcal{J}^{-\alpha}_{t,b}f(.), \alpha \in \mathbb{R}^{-}/\mathbb{Z}^{-}$ is the right-sided Riemann-Liouville fractional integral of order $-\alpha$ defined in \cite{podlubny2007adjoint} by:}
\R{\begin{equation*}
\mathcal{J}^{\beta}_{t,b}f(t)=\frac{(-1)^{n}}{\Gamma(\beta)} \int _{t} ^{b} (t-\tau)^{\beta-1} f(\tau)d\tau,   
\\
\end{equation*}}
\R{where, $n-1 < \beta \leq n, t\in [a,b], f(.) \in C[a,b]$}
\end{definition}
\begin{definition}\label{fourier} 
Fourier transform definition, its inverse and the Fourier transform of the classical derivative are given by:\\
\R{Let $f(t) \in C^{m}(\mathbb{R}), \mathcal{F}  \{f(t) \}(s) \in L^{1}(\mathbb{R})$ then, we have:}
\begin{equation*}
\left\{\begin{array}{llll}    
\mathcal{F}  \{f(t) \}(s)=\int_{-\infty}^{+\infty} e^{ist}f(t)dt, \\
\mathcal{F}^{-1}  \{\tilde{f}(s) \}(t)=\frac {1}{2\pi}\int_{-\infty}^{+\infty} e^{-ist}\tilde{f}(s)ds, \\
 \mathcal{F}  \{D^{m}_{t}f(t) \}(s)=(-is)^{m} \mathcal{F}  \{f(t) \}(s),
\end{array}\right.
\end{equation*}
 where, $\mathcal{F}  \{f(t) \}(s):=\tilde{f}(s), \R{s \in \mathbb{R}}, m \in {\rm I\!N}, i^{2}=-1$. 
\end{definition}
\begin{lemma}
 Fourier transform of Caputo fractional derivative \R{is} given in \cite{huang2005fundamental} by: \R{Let $f(t) \in C^{n}(\mathbb{R})$, then:}
\begin{equation} \label{eq1} 
\mathcal{F}  \{^{C}D^{\alpha}_{t}f(t) \}(s)=(-is)^{\alpha} \mathcal{F}  \{f(t) \}(s) 
\end{equation}
 where, $n-1< \alpha < n, n \in {\rm I\!N}$.
\end{lemma}

\begin{lemma}
Fractional integration by parts was given in  \cite{podlubny2007adjoint}  by:\\
 \R{Let $f(t),g(t) \in C^{n}(\mathbb{R})$, then:}
\vspace{-1.8cm}
\begin{equation}\label{eq2} 
\begin{split}
\int_{a}^{b} f(t)^{C}_{a}D^{\alpha}_{t}g(t)dt \qquad \qquad \\
=\sum _{k=0}^{n-1} (-1)^{n-1-k} [g^{(k)}(t)\tensor*[^{C}]{D}{^{\alpha-1-k}_{t,b}}f(t)]_{t=a}^{t=b}\\
+(-1)^{(n)} \int _{a}^{b} g(t) \tensor*[^{C}]{D}{^{\alpha}_{t,b}}f(t) dt 
\end{split}
\end{equation}

where, $n-1< \alpha < n, n \in {\rm I\!N}$. $\tensor*[^{C}]{D}{^{\alpha}_{t,b}}$ is the right-sided Caputo derivative defined in \eqref{sided}.
\end{lemma}
\vspace{-0.3cm}

\begin{theorem}\label{thm20}  {\cite{podlubny1998fractional} }
Let $t\in [a,b],$ we have: 
\begin{equation}\label{gama}
^{C}_{\R{t_{0}}}D^{\alpha}_{t} (t-t_{0})^{\beta}=\frac{\Gamma (\beta+1)}{\Gamma (\beta-\alpha+1)}(t-t_{0})^{\beta-\alpha}, \quad \forall \beta 
\end{equation} 
where $t_{0} \in [a,b], n-1 < \alpha \leq n, n \in {\rm I\!N},$ see \cite{podlubny1998fractional}  for the proof.
\end{theorem}

\begin {theorem}\label{eq3}
Caputo Fractional derivative of an integral depending on a parameter is given in  \cite{ghaffour_laleg_kirati}  by:\\
\R{Let  $K(t,\tau)$ be a function such that both $K(t,\tau)$ and its partial derivative $\frac{\partial}{\partial t}K(t,\tau)$ are continuous in $t$ and $\tau$ then, we have:}
\vspace{-0.3cm}
\begin{equation*}
^{C}D^{\alpha}_{t}\int_{0}^{t}K(t, \tau)d\tau=\int_{0}^{t} \hspace{1mm}  _{\tau}^{C}D^{\alpha}_{t}K(t, \tau)d\tau+\R{\mathcal{J}^{1-\alpha}_{t}K(t, t)}, 
\end{equation*}
where, $0 \leq \alpha < 1$, \R{$K(t,\tau)$ is defined over $\mathbb{R}^{+}_{*}X\mathbb{R}^{+}_{*}$  \text{such that} $0<\tau<t$.}
\end {theorem}
\begin {proof}
Using the \R{Leibniz} differentiation rule for the classical first derivative, similar reasoning to the proof in \cite{podlubny1998fractional} gives:
\begin{equation*} 
\begin{split}
 ^{C}D^{\alpha}_{t} \int_{0}^{t}K(t, \tau)d\tau=\frac{1}{\Gamma(1-\alpha)} \int_{0}^{t}  \frac{1}{(t-s)^{\alpha}}  \frac{\partial}{\partial s} \int _{0}^{s} K(s,\tau)d\tau ds \\
=\frac{1}{\Gamma(1-\alpha)} \int_{0}^{t}  \frac{1}{(t-s)^{\alpha}} \left[  \int _{0}^{s} \frac{\partial}{\partial s} K(s,\tau)d\tau+ K(s, s) \right] ds\\
 =\frac{1}{\Gamma(1-\alpha)} \int_{0}^{t}  \frac{1}{(t-s)^{\alpha}} \int _{0}^{s} \frac{\partial}{\partial s} K(s,\tau)d\tau ds+\frac{1}{\Gamma(1-\alpha)}  \int_{0}^{t}  \frac{1}{(t-s)^{\alpha}} K(s, s) ds  \\
=\frac{1}{\Gamma(1-\alpha)} \int_{0}^{t} \int _{\tau}^{t} \frac{1}{(t-s)^{\alpha}}  \frac{\partial K(s,\tau)}{\partial s}  ds  d\tau+ \hspace{-2mm} _{\quad}\hspace{-1mm}^{C}D^{\alpha-1}_{t}K(t,t) \\
=\int_{0}^{t} \hspace{1mm } _{\tau}^{C}D^{\alpha}_{t}K(t, \tau)d\tau+ \mathcal{J}^{1-\alpha}_{t}K(t,t).
\end{split}
\end{equation*} 
which completes the proof.
\end {proof}

\begin{definition}\label{AS}{\cite{beck2012brief}}
\R{
Consider the PDE $ \frac{\partial}{\partial t}u=Lu , u \in X$ with $L$ is some linear operator on the Banach space $X$.\\
A solution u(t) of this system is said to be stable if, given any $\epsilon > 0,$ there exists a $\delta=\delta(\epsilon) >0$  such that, for all $u_{0} \in X$ with $\left\|u_{0}\right\|_{X} \leq \delta$, the corresponding solution satisfies $\left\|u(t)\right\|_{X} \leq \epsilon$ for all $t \geq 0$. If in addition there exists a $\delta^{*}$
such that for all initial conditions with $\left\|u_{0}\right\|_{X} \leq \delta^{*}$
the corresponding solution satisfies \[ \lim_{t \rightarrow \infty} \left\|u(t)\right\|_{X} =0, \] then this  solution is said to be asymptotically
stable.}
\end{definition}
\begin{remark}
\R{In what follows, we will be considering the stability and the asymptotic stability in $L_{\infty}(\mathbb{R})$. That is, taking $X=L_{\infty}(\mathbb{R})$ in definition \ref{AS}.}
\end{remark}
\begin{remark}
\R{In what follows, We will be using the following norms: 
$$
\norm{P(x,t)}_{\infty}=\sup_{x \in [0,L]} |P(x,t)|
, \quad \quad \text{and } \quad  \norm{P(x,t)}_{2}=(\int_{0}^{L}|P(x,t)|^{2} dx)^{1/2}.
$$}
\R{which are defined over $L^{\infty}(\mathbb{R})$ and $L^{2}(\mathbb{R})$ respectively. } 
\end{remark}

\begin{theorem}\label{th2.11}
\R{
Let $f(t) \in L^{1}(\mathbb{R})$. Then the Fourier transform of $f(t)$ is bounded in $L^{1}(\mathbb{R})$ by:
$$
|\mathcal{F}  \{f(t) \}(s)|\leq \int_{-\infty}^{+\infty} |f(t)|dt=\norm{f(t)}_{L^{1}(\mathbb{R})}.
$$}
\end{theorem}
\begin{theorem}\label{th2.12}
\R{Fourier transform of  fractional integral is given in \cite{podlubny1998fractional} by:
\begin{equation*}
\mathcal{F}  \{\mathcal{J}^{\alpha}_{t}f(t) \}(s)=\frac{1}{(is)^{\alpha}} \mathcal{F}  \{f(t) \}(s) , \quad s\neq0\\ 
\end{equation*}
 where, $n-1< \alpha < n, n \in {\rm I\!N}$.\\}

\end{theorem}
 \section{Problem Statement}

\begin{figure}[ht]
	\centering
         \includegraphics[scale=0.29]{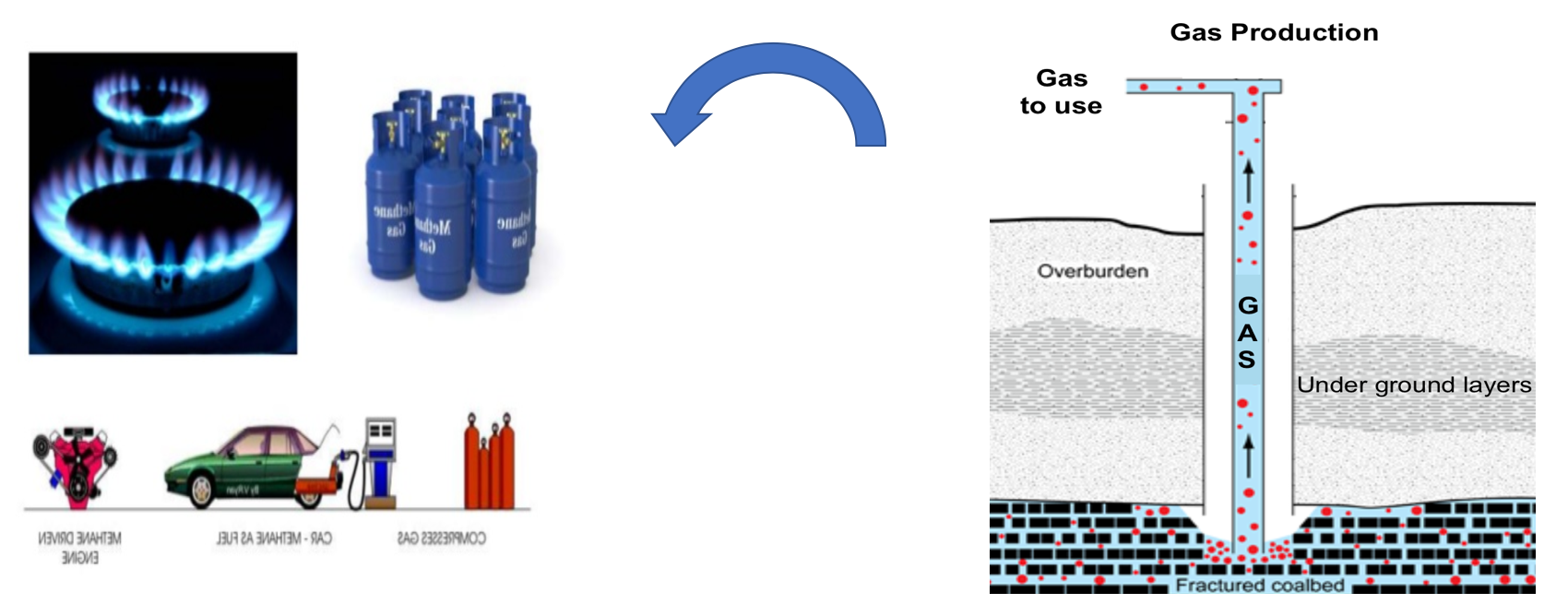}
	\caption{Gas production process}\label{pic1}    
\end{figure}

This section presents the model's equations. We consider the gas production mechanism which is done by drilling wells. When a vertical well is drilled through a gas saturated rock, the underground layers burden pressure causes the gas to flow through the drilled well into the production pipelines as shown in figure \ref{pic1}. Table \ref{Tab1} summarizes the system's parameters.
\begin{table}[h]
	\caption{Parameters Description}
	\label{Tab1}
	\begin{center}
		\begin{tabular}{|c||c|}
			\hline
			Symbol  & Description \\
			\hline
			$t$  & Time \\ 
			$x$  & Space \\
			$P(x, t)$ & Gas pressure \\
			$\rho(x,t)$ &  Density \\
			$\varphi(x,t)$ &  Porosity \\
			$u(x,t)$ &  Velocity \\
			$k(t)$ &  Permeability  \\
			$\mu$ &  Viscosity  \\
			$\alpha$ & Diffusion coefficient \\
			$Q(x,t)$ & Gas production flow  \\
		         $C$ & Variation of the porosity \\   
		        $ $ &  with respect to the pressure\\
		        $\R{u(t)}$ & \R{Control input} \\
			\hline
		\end{tabular}
	\end{center}
\end{table}

\subsection{Gas Production Fractional Model}
%We start by recalling the classical PDE derivation in \cite{amir2018physics} for the gas pressure :
%\begin{enumerate}
%\item Mass conservation law :  
%           \begin{equation*}
%         \frac{\partial (\rho \varphi)(P(x,t))}{\partial t} +\nabla. (\rho u)(P(x,t))=Q(x,t),
%           \end{equation*}

%\item Assume the non-compressible fluid, which means that the density does not change with respect to the pressure:
 %          \begin{equation*}
  %       \rho \frac{\partial \varphi(P(x,t))}{\partial P(x,t)}  \frac{\partial P(x,t)}{\partial t}  +\rho \nabla. u(P(x,t))=Q(x,t),
   %        \end{equation*}
%\item Assume a constant-compressible rock (a particular case of the linear-compressibility), which means that the porosity is constant with respect to the pressure:
 %          \begin{equation*}
  %       \frac{\partial P(x,t)}{\partial t}  +\frac{1}{C} \nabla. u(P(x,t))=\frac{1}{C \rho}Q(x,t), 
   %        \end{equation*}
%\item Using the Darcy flow equation:           
 %         \begin{equation*}
  %       \frac{\partial P(x,t)}{\partial t}  +\frac{1}{C} \nabla. (-\frac{k(t)}{\mu} \frac{\partial P(x,t)}{\partial x})=\frac{1}{C \rho}Q(x,t).
   %        \end{equation*}
%\end{enumerate}
The standard modeling approach for gas pressure is recalled in the appendix \R{A}.  However, due to the rocks explosion caused by the gas pressure, the constructed pores don't have the same shapes and dimensions. Thus, the variation of the pressure with respect to space is not equally distributed which means that the Darcy law equation doesn't describe fully the considered phenomena. Therefore,  the authors in \cite{amir2018physics} suggested to use the non-Darcy law, where the first spatial derivative of the pressure is replaced with a fractional spatial derivative of order $\alpha$, where, $0< \alpha \leq 1$ represents the diffusion coefficient. Therefore the FPDE that fully describes the considered phenomena is given with the boundary conditions and the system control by:       
\begin{equation}\label{eq4} 
\left\{\begin{array}{llll}    
 \dfrac{\partial P(x,t)}{\partial t}  \R{-}\dfrac{1}{C}  \dfrac{k(t)}{\mu} \dfrac{\partial}{\partial x} \hspace{-3mm}\quad^{\R{C}}D^{\alpha}_{x} P(x,t)=\dfrac{1}{C \rho}Q(x,t),\\
\\
P(x, 0) = g_{0}(x),     \\ 
P_{x}(0, t) = 0,     P(L, t) = u(t), \\

\end{array}\right.
\end{equation}
where $t > 0, x \in [0,L]$, $u(t)$ is the control input, $P(x, t)$ is the gas pressure distributed in space and in time. The permeability $k(t)$ and the viscosity $\mu$ are both positive. $\tensor*[^{C}]{D}{^{\alpha}_{x}} P(x,t)$ is the Caputo space fractional derivative of order $\alpha$. \R{$u(t)$ is the control input.}

\begin{remark} 
It has been proved in \cite{podlubny1998fractional} that if we consider the Caputo fractional derivative then the operator
 $\frac{\partial}{\partial x} \tensor*[^{C}]{D}{^{\alpha}_{x}} P(x,t) $in \eqref{eq4} cannot be written as $ \tensor*[^{C}]{D}{^{\beta}_{x}} P(x,t)$ with $1< \beta \leq 2$.
\end{remark} 

The main objective of this paper is to track the pressure gradient at final position in the presence of disturbances using some measurements and a boundary control. However, before that, we first solve the problem analytically. Then, we study the well-posedness and the stability of the considered system.

\section {Fundamental Solution of the Space FPDE}
In this section, we first derive the fundamental solution for the non homogeneous Caputo space fractional PDE. For this purpose, we use the Fourier transform of Caputo fractional derivative  \eqref{eq1}. Then, we investigate the uniqueness of the obtained solution.

\subsection{Analytic Solution}
\begin{theorem} \label{thm1}
\R{Assuming that  $Q(x,.),g_{0}(x)\in \mathcal{L}^{1} (\mathbb{R})$ and $u(t) \in C^{1}(\mathbb{R^{+}})$,}, system \eqref{eq4} admits at least one solution and it is given by:
\begin{equation*} 
\begin{split}
P(x, t)=u(t)+ \int _{- \infty}^{+ \infty}  \int _{0}^{t} \mathcal{G}_{\alpha}(\mid x-y \mid,t-\tau)  \bar{Q}(y,\tau) d\tau dy 
+ \int _{- \infty}^{+ \infty}  \mathcal{G}_{\alpha}(\mid x-y \mid,t) \bar{g}_{0}(y) dy. 
\end{split}
\end{equation*}
where, 
$
 \bar{g}_{0}(x):=g_{0}(x)-u(0), \bar{Q}(x,t):=\frac{1}{C \rho} Q(x,t)-\frac{\partial}{\partial t}u(t)
$ \R{ and 
$\mathcal{G}_{\alpha}(x,t)$ is the Green function defined by:
$$
\mathcal{G}_{\alpha}(x,t)=\frac{1}{2 \pi}\int_{-\infty}^{+\infty} e^{-isx}e^{\frac{1}{C}\frac{k(t)}{\mu} \phi^{\alpha+1}(s)t} ds,
$$
with $\phi^{\alpha+1}(s)=(-is)^{\alpha+1}. $ }

\end{theorem}
\begin{proof}
Starting from the system \eqref{eq4}, we define the new coordinate $\bar{P}(x,t)$,
\begin{equation}\label{BCs}
\bar{P}(x,t):=P(x,t)-u(t),
\end{equation}
where, $u(t)$ is the boundary condition in system \eqref{eq4}. The system for the new coordinate becomes:\\
\begin{equation}\label{e5}
\left\{\begin{array}{llll}    
 \dfrac{\partial \bar{P}(x,t)}{\partial t}  \R{-}\dfrac{1}{C}  \dfrac{k(t)}{\mu} \dfrac{\partial}{\partial x} \tensor*[^{C}]{D}{^{\alpha}_{x}} \bar{P}(x,t)=\bar{Q}(x,t),\\
\\
\bar{P}(x, 0) = \bar{g}_{0}(x),     \\ 
\bar{P}_{x}(0, t) = 0,     \bar{P}(L, t) = 0,     \\

\end{array}\right.
\end{equation}
where, $\bar{g}_{0}(x)$ is the new initial condition and $\bar{Q}(x,t)$ is the new source term given respectively by: 
$$
 \bar{g}_{0}(x):=g_{0}(x)-u(0) \quad  \bar{Q}(x,t):=\frac{1}{C \rho}Q(x,t)-\frac{\partial}{\partial t}u(t).
$$
Then, we multiply the gas pressure $\bar{P}(x,t)$ in \eqref{eq4} by the indicator function of $[0,L]$ to have the FPDE defined for $x \in  \mathbb{R}$.  Then we apply the superposition principle to write the solution of \eqref{eq4} as $\bar{P}=h+v$, where, for every $t > 0, x \in \mathbb{R}$, $u$ and $v$ are the solution of the following problems respectively:
\begin{equation}\label{eq5} 
\left\{\begin{array}{llll}    
 \dfrac{\partial h(x,t)}{\partial t}  \R{-}\frac{1}{C}  \dfrac{k(t)}{\mu} \dfrac{\partial}{\partial x}\tensor*[^{C}]{D}{^{\alpha}_{x}} h(x,t)=0,  \\
\\
h(x, 0) = \bar{g}_{0}(x),     \\ 
\end{array}\right.
\end{equation} \\
\vspace{-4mm}
\begin{equation}\label{eq6} 
\left\{\begin{array}{llll}    
 \dfrac{\partial v(x,t)}{\partial t}  \R{-}\frac{1}{C}  \dfrac{k(t)}{\mu} \dfrac{\partial}{\partial x}_{}^{C}D^{\alpha}_{x} v(x,t)=\bar{Q}(x,t) \\
\\
v(x, 0) = 0,     \\ 
\end{array}\right.
\end{equation}

\R{because, both the Caputo fractional derivative and classical derivative are linear.} We start by solving system \eqref{eq6}, we apply the Duhamel's principle \cite{jeffrey2003applied} \R{( which allows to solve first the homogeneous linear, PDE, and then superposing to find the solution of he original PDE)}, the solution of \eqref{eq6} is given by:

\begin{equation}\label{eq7} 
v(x, t)=\int _{0}^{t} \mathbb{V}(x,t,\tau)d\tau,  
\end{equation}

where, for all $\tau \in [0,t]$, for every $t > 0, x \in \mathbb{R}$. \R{We start by computing the time first derivative of $v(x,t)$ in \eqref{eq7}:
\begin{equation}\label{V*} 
\begin{split}
\frac{\partial}{\partial t} v(x, t)=\frac{\partial}{\partial t} \int _{0}^{t} \mathbb{V}(x,t,\tau)d\tau,  \\
= \int _{0}^{t} \frac{\partial}{\partial t} \mathbb{V}(x,t,\tau)d\tau +\mathbb{V}(x,t=\tau),  \\
\end{split}
\end{equation}
using the classical Leibnitz differentiation rule. Now, we compute the space fractional derivative of order $\alpha$ of $v(x,t)$ in \eqref{eq7}:
\begin{equation}\label{V**} 
\begin{split}
^{C}D^{\alpha}_{x} v(x, t)=^{C}D^{\alpha}_{x} \int _{0}^{t} \mathbb{V}(x,t,\tau)d\tau,  \\
= \int _{0}^{t} {}^{C}D^{\alpha}_{x} \mathbb{V}(x,t,\tau)d\tau  \\
\Rightarrow \frac{\partial}{\partial x}^{C}D^{\alpha}_{x} v(x, t)= \int _{0}^{t}\frac{\partial}{\partial x}^{C}D^{\alpha}_{x} \mathbb{V}(x,t,\tau)d\tau 
\end{split}
\end{equation}
Thus, }$\mathbb{V}(.,.,\tau)$ satisfies the following system:

\begin{equation}\label{eq8} 
\left\{\begin{array}{llll}    
\frac{\partial \mathbb{V}(x, t)}{\partial t}  \R{-}\frac{1}{C}  \frac{k(t)}{\mu} \frac{\partial}{\partial x}_{}^{C}D^{\alpha}_{x}\mathbb{V}(x, t)=0,\\
\\
\lim_{t \to \tau} \quad  \mathbb{V}(x, t)= \bar{Q}(x,\tau).     \\ 
\end{array}\right.
\end{equation}
Then, by applying the Fourier transform to the FPDE with respect to the variable $x$, we have:
\begin{equation*} 
\left\{\begin{array}{llll} 
 \mathcal{F} \{ \frac{\partial \mathbb{V}(x, t)}{\partial t} \} \R{-}\frac{1}{C}  \frac{k(t)}{\mu}  \mathcal{F}\{ \frac{\partial}{\partial x} _{}^{C}D^{\alpha}_{x} \mathbb{V}(x, t)\} =0, \\
 \R{\mathcal{F} \{ \mathbb{V}(x, t=\tau) \}= \mathcal{F} \{ \bar{Q}(x,\tau)\},}\\
 \end{array}\right.
\end{equation*}
using the Fourier transform of the first spatial derivative \R{(definition \ref{fourier})} and Caputo fractional derivative \eqref{eq1}, we get:
\begin{equation}\label{8}  
\left\{\begin{array}{llll} 
\frac{\partial \tilde{\mathbb{V}}(s, t)}{\partial t}  \R{-}\frac{1}{C}  \frac{k(t)}{\mu} (-is)^{\alpha+1} \tilde{\mathbb{V}}(s, t)=0, \\
\R{\tilde{ \mathbb{V}}(s, t=\tau) =   \tilde{\bar{Q}}(x,\tau),}\\
 \end{array}\right.
\end{equation}
\R{where, $\tilde{\bar{Q}}(s,\tau)=\mathcal{F} \{ \bar{Q}(x,\tau) \}(s,\tau)$.} Using the initial condition in \eqref{eq8}, the solution of \eqref{8} is given by:
\begin{equation} \label{9}  
\tilde{\mathbb{V}}(s, t)=e^{  \frac{1}{C} \frac{k(t)}{\mu} (-is)^{\alpha+1} (t-\tau)}\tilde{\bar{Q}}(s,\tau), \\
\end{equation}

where, $\tilde{\bar{Q}}(s,\tau)=\mathcal{F} \{ \bar{Q}(x,\tau) \}(s,\tau)$. Applying the inverse Fourier transform \eqref{fourier}, \eqref{9} becomes:
\begin{equation*} 
\mathbb{V}(x, t)=
\frac{1}{2\pi}\int_{-\infty}^{+\infty} e^{-isx}e^{  \frac{1}{C} \frac{k(t)}{\mu} (-is)^{\alpha+1} (t-\tau)}\tilde{\bar{Q}}(s,\tau)ds,
\end{equation*}
which can be written as:
\begin{equation} \label{eq9}  
\mathbb{V}(x, t)=\frac{1}{2\pi}\int _{- \infty}^{+ \infty} e^{-isx}  \tilde{\mathcal{G}_{\alpha}}(s,t-\tau)\tilde{\bar{Q}}(s,\tau)ds, \\
\end{equation}
where,  $\mathcal{G}_{\alpha}(x,t)$ is the Green function defined in \cite{huang2005fundamental} by:
$$
\mathcal{G}_{\alpha}(x,t)=\frac{1}{2 \pi}\int_{-\infty}^{+\infty} e^{-i\R{s}x}e^{\frac{1}{C}\frac{k(t)}{\mu} \phi^{\alpha+1}(s)t} d\R{s}
$$
and characterized by: 
\begin{equation} \label{eq10}  
\tilde{\mathcal{G}_{\alpha}}=\mathcal{F}\{  \mathcal{G}_{\alpha} \} (s,t)=e^{\frac{1}{C}\frac{k(t)}{\mu} \phi^{\alpha+1}(s)t} ,
\end{equation}
where,
$$\phi^{\alpha+1}(s)=(-is)^{\alpha+1}.  $$
$\mathcal{F}$ is the Fourier Transform, more details in the green function can be found in \cite{aldoghaither2017direct}, \cite{huang2005fundamental}. Then, using the inverse Fourier transform \eqref{fourier}, \eqref{eq9} becomes:
\begin{equation*}   
\mathbb{V}(x, t)=\int _{- \infty}^{+ \infty}   \mathcal{G}_{\alpha}( x-y ,t-\tau) \bar{Q}(y,\tau) dy , \\
\end{equation*}

Thus, using \eqref{eq7}, the solution of \eqref{eq6} is given by:
\begin{equation} \label{eq11}  
v(x, t)=\int _{- \infty}^{+ \infty}  \int _{0}^{t} \mathcal{G}_{\alpha}( x-y ,t-\tau)\bar{Q}(y,\tau)) d\tau dy \\
\end{equation}
Similarly, the solution of \eqref{eq5} is given by:
\begin{equation} \label{eq12}  
h(x, t)=\int _{- \infty}^{+ \infty}  \mathcal{G}_{\alpha}( x-y 
,t)\bar{g}_{0}(y) dy. \\
\end{equation}
Finally, the fundamental solution of \eqref{eq4} is:
\begin{equation} \label{eq13}  
\begin{split}
P(x, t)=u(t)+\int _{- \infty}^{+ \infty}  \int _{0}^{t} \mathcal{G}_{\alpha}( x-y ,t-\tau) \bar{Q}(y,\tau) d\tau dy 
+ \int _{- \infty}^{+ \infty}  \mathcal{G}_{\alpha}( x-y ,t)\bar{g}_{0}(y) dy.  
\end{split}
\end{equation}
\end {proof}
\subsection{Uniqueness}
The uniqueness of the solution is a direct result from the linearity of the integral operator.
\begin{theorem} \label{thm22}
\R{Assuming that  $g_{0}(x)\in \mathcal{L}^{1} (\mathbb{R})$ and $u(t) \in C^{1}(\mathbb{R^{+}})$, Let $P_{1}(x, t)$ and $P_{2}(x, t)$ be two solutions of system \eqref{eq4} with different source terms $Q_{1}(x,.),Q_{2}(x,.) \in \mathcal{L}^{1} (\mathbb{R})$ respectively. Then, the
condition $Q_{1}(x,t)=Q_{2}(x,t)$ implies that $P_{1}(x, t)=P_{2}(x, t)$.}
\end{theorem}
\begin{proof}
\R{ Suppose that system \eqref{eq4} admits two solutions $P_{1}(x, t)$ and $P_{2}(x, t)$ with different source terms $Q_{1}(x,.),Q_{2}(x,.) \in \mathcal{L}^{1} (\mathbb{R})$. Then, by Theorem \ref{thm1}, $P_{1}(x, t)$ and $P_{2}(x, t)$ are both given by:
\begin{equation} \label{p1}  
\begin{split}
P_{1}(x, t)=\int _{- \infty}^{+ \infty}  \hspace{-2mm}\int _{0}^{t}\hspace{-0.4mm} \mathcal{G}_{\alpha}( x-y ,t-\tau) Q_{1}(y,\tau) d\tau dy
& + \int _{- \infty}^{+ \infty} \hspace{-4mm} \mathcal{G}_{\alpha}(\mid x-y \mid,t)\bar{g}_{0}(y) dy+u(t).  \\
\end{split}
\end{equation}
and 
\begin{equation} \label{p2}  
\begin{split}
P_{2}(x, t)=\int _{- \infty}^{+ \infty}  \int _{0}^{t} \mathcal{G}_{\alpha}( x-y ,t-\tau) Q_{2}(y,\tau)  d\tau dy 
&+ \int _{- \infty}^{+ \infty}  \mathcal{G}_{\alpha}( x-y ,t)\bar{g}_{0}(y) dy+u(t).  \\
\end{split}
\end{equation}
Then, by \eqref{p1} and \eqref{p2}:
\begin{equation*}  
P_{1}(x, t)-P_{2}(x, t)=\int _{- \infty}^{+ \infty}  \int _{0}^{t} \mathcal{G}_{\alpha}( x-y ,t-\tau) \big[Q_{1}(y,\tau)-Q_{2}(y,\tau)\big] d\tau dy  
\end{equation*}
completes the proof.}
\end{proof}

\section{Stability Conditions}
In this section, we aim to derive some sufficient conditions to guarantee the stability of the problem given in \eqref{eq4}. 
\begin{theorem}\label{thm2} 
\R{\hspace{-2mm}Suppose that $\bar{Q}(x,t)$  is separable on its variables $\bar{Q}(x,t)\hspace{-0.8mm}=T(t)q(x)$ such that $q(x) \in L^{1}(\mathbb{R})$ and $T(t)$ is bounded on time, this  means that 
$$
\bar{Q}(x,t) \leq r q(x)
.$$ 
Suppose as well that $\bar{g}_{0}(x) \in L^{1}(\mathbb{R})$ and $k(t)$ is bounded from bellow $(k_{0} \leq k(t))$ and that the control $u(t)$ is bounded as well.  Then, the solution of system \eqref{eq4} is stable in $\B{L_{\infty}(\mathbb{R})}$. That is:
$$
\lim_{t \to +\infty} \quad \norm{ P(x, t)}_{\infty} < \infty
$$\\}
\end{theorem}
\vspace{-2mm}
\begin{proof}
\R{
Using \R{T}heorem \ref{thm1}, the solution of \eqref{eq4} is given by \eqref{eq13}. $\bar{Q}(x,t)$ and $u(t)$ are bounded in time which means that:
\begin{equation*} 
 \bar{Q}(x,t) \leq r q(x).
\end{equation*}
\begin{equation*} 
u(t)\leq u_{max}.
\end{equation*}
 Thus, \eqref{eq4} becomes:
 \begin{equation*} 
\begin{split}
P(x, t)= u(t) + \int _{- \infty}^{+ \infty}  \int _{0}^{t} \mathcal{G}_{\alpha}( x-y ,t-\tau)   \bar{Q}(y,t)  d\tau dy 
+  \int _{- \infty}^{+ \infty}  \mathcal{G}_{\alpha}( x-y ,t)\bar{g}_{0}(y) dy .  \\
= \frac{1}{2 \pi} \int _{-\infty}^{+\infty} \int _{0}^{t}   \tilde{\mathcal{G}}_{\alpha}(s ,t-\tau)\tilde{\bar{Q}}(s,t)  e^{-isx}  d\tau ds 
+\frac{1}{2 \pi} \int _{-\infty}^{+\infty}  \tilde{\mathcal{G}}_{\alpha}(s,t)\tilde{\bar{g}}_{0}(s) e^{-isx}ds+u(t),  \\
\end{split}
\end{equation*}
because both $\bar{Q}(x,t)$ and  $\bar{g}_{0}(x)$ are in $L^{1}(\mathbb{R})$ which means that their Fourier transforms $\tilde{\bar{Q}}(x,t)$ and  $\tilde{\bar{g}}_{0}(x)$ exist. Then, we get:
\begin{equation} \label{eq17}  
\begin{split}
\big|P(x, t)\big| &   \leq \frac{1}{2 \pi} \big| \int _{-\infty}^{+\infty}  \int _{0}^{t}  \tilde{\mathcal{G}}_{\alpha}( s ,t-\tau) e^{-isx}\tilde{\bar{Q}}(s,t)    d\tau   ds \big|
+ \frac{1}{2 \pi} \big| \int _{-\infty}^{+\infty}     \tilde{\mathcal{G}}_{\alpha}(s,t) e^{-isx}\tilde{\bar{g}}_{0}(s)    ds\big|+\big|u_{max}|,\\
&\leq\frac{1}{2 \pi} \int _{-\infty}^{+\infty} \big| \int _{0}^{t}  \tilde{\mathcal{G}}_{\alpha}( s ,t-\tau) e^{-isx}\tilde{\bar{Q}}(s,t)    d\tau  \big|  ds
+ \frac{1}{2 \pi} \int _{-\infty}^{+\infty}  \big|   \tilde{\mathcal{G}}_{\alpha}(s,t) e^{-isx}\tilde{\bar{g}}_{0}(s)  \big|  ds+\big|u_{max}|,\\
&\leq \frac{1}{2\pi} \int _{-\infty}^{+\infty} \big| r \tilde{q}(s)e^{-isx} \big| \big| \int _{0}^{t} \tilde{\mathcal{G}}_{\alpha}(s,t-\tau) d\tau\big|  ds 
+ \frac{1}{2\pi}\int _{-\infty}^{+\infty}  \big| \tilde{\mathcal{G}}_{\alpha}(s,t)\big| \big|\tilde{\bar{g}}_{0}(s) \big|ds+\big|u_{max}\big|.  \\
\end{split}
\end{equation}
By taking the sup on $x$ over $[0,L]$ of \eqref{eq17}, we get:
\begin{equation} \label{eq18}  
\begin{split}
 \norm{ P(x, t) }_{\infty} \leq 
  \frac{|r|}{2 \pi} \int _{-\infty}^{+\infty} \big|\tilde{q}(s)\big|  \big| \int _{0}^{t}  \tilde{\mathcal{G}}_{\alpha}( s ,t-\tau) d\tau\big|   ds 
+\frac{1}{2 \pi}  \int _{-\infty}^{+\infty}  \big|\tilde{\mathcal{G}}_{\alpha}(s,t)\big| \big| \tilde{\bar{g}}_{0}(s) \big| ds+\big|u_{max}\big|,  \\
\end{split}
\end{equation}
taking the limit of \eqref{eq18}, we get:
\begin{equation} \label{eq19}  
\begin{split}
\lim_{t \rightarrow +\infty} \norm{ P(x, t) }_{\infty} \leq 
  \frac{|r|}{2 \pi} \int _{-\infty}^{+\infty} \big| \tilde{q}(s)\big| \lim_{t \rightarrow +\infty}  \big| \int _{0}^{t}  \tilde{\mathcal{G}}_{\alpha}(s ,t-\tau) d\tau\big|   ds +\frac{1}{2 \pi} \int _{-\infty}^{+\infty}   \lim_{t \rightarrow +\infty}   \big|\tilde{\mathcal{G}}_{\alpha}(s,t)\big| \big| \tilde{\bar{g}}_{0}(s) \big| ds+\big|u_{max}\big|,  \\
\end{split}
\end{equation}
where
\begin{equation} \label{eq10*}  
\tilde{\mathcal{G}}_{\alpha}(s,t)= e^{\frac{1}{C}\frac{k(t)}{\mu} \phi^{\alpha+1}(s)t}
\end{equation}
and
$$\phi^{\alpha+1}(s)=(-is)^{\alpha+1}.  $$\\
On the other hand, we have that, for all $0<\alpha \leq 1$: 
\begin{equation} \label{ao}  
\begin{split}
\phi^{\alpha+1}(s)=(-is)^{\alpha+1}=s^{\alpha+1}(-i)^{\alpha+1}=s^{\alpha+1}e^{-i\frac{\pi}{2}(\alpha+1)}&\\
=s^{\alpha+1} \{ \cos(\frac{\pi}{2}(\alpha+1))-i \sin(\frac{\pi}{2}(\alpha+1)) \} &\\
=s^{\alpha+1} \{ -\sin(\frac{\pi}{2} \alpha)-i \cos(\frac{\pi}{2} \alpha) \}. &\\
\Rightarrow \operatorname{Re}\{ \phi^{\alpha+1}(s) \} \leq 0, \quad   s\geq 0.
\end{split}
\end{equation}
if $s <0$ we have that, for all $0<\alpha \leq 1$: 
\begin{equation} \label{aoo}  
\begin{split}
\phi^{\alpha+1}(s)=(-is)^{\alpha+1}=s^{\alpha+1}(-i)^{\alpha+1}=s^{\alpha+1}e^{-i\frac{\pi}{2}(\alpha+1)}&\\
=s^{\alpha+1} \{ \cos(\frac{\pi}{2}(\alpha+1))-i \sin(\frac{\pi}{2}(\alpha+1)) \} &\\
=s^{\alpha+1} \{ -\sin(\frac{\pi}{2} \alpha)-i \cos(\frac{\pi}{2} \alpha) \}. &\\
=-(-1)^{\alpha}(-s)^{\alpha+1} \{ -\sin(\frac{\pi}{2} \alpha)-i \cos(\frac{\pi}{2} \alpha) \}. &\\
=-(i)^{2\alpha}(-s)^{\alpha+1} \{ -\sin(\frac{\pi}{2} \alpha)-i \cos(\frac{\pi}{2} \alpha) \}. &\\
=-e^{i \pi \alpha}(-s)^{\alpha+1} \{ -\sin(\frac{\pi}{2} \alpha)-i \cos(\frac{\pi}{2} \alpha) \}. &\\
=(-s)^{\alpha+1} \{ -\cos(\pi \alpha)-i \sin(\pi \alpha) \} \{ -\sin(\frac{\pi}{2} \alpha)-i \cos(\frac{\pi}{2} \alpha) \}. &\\
=(-s)^{\alpha+1} \{ \cos(\pi \alpha)sin(\frac{\pi}{2} \alpha)- \sin(\pi \alpha)\cos(\frac{\pi}{2} \alpha)\\
+i (\cos(\frac{\pi}{2} \alpha)\cos(\pi \alpha)+\sin(\frac{\pi}{2} \alpha)\sin(\pi \alpha))   \}. &\\
=(-s)^{\alpha+1} \{- \sin  (\frac{\pi}{2} \alpha)+i \cos  (\frac{\pi}{2} \alpha) \}. &\\
\Rightarrow \operatorname{Re}\{ \phi^{\alpha+1}(s) \} \leq 0, \quad   s< 0.
\end{split}
\end{equation}
Thus from \eqref{ao} and \eqref{aoo} we have that:
\begin{equation} \label{a}  
\operatorname{Re}\{ \phi^{\alpha+1}(s) \} \leq 0, \quad 0<\alpha \leq 1, \forall s \in \mathbb{R}, 
\end{equation}
using the fact that both the permeability and the viscosity are both positive and the linear-compressibility with $C > 0$, the sign of the real part of $\phi^{\alpha+1}(s)$ is negative. Using \eqref{eq10*} we have that:
\begin{equation} \label{eq20}  
\begin{split}
\big|\tilde{\mathcal{G}}_{\alpha}( s ,t)\big| =   e^{\frac{1}{C}\frac{k(t)}{\mu} \operatorname{Re} \{ \phi^{\alpha+1}(s) \}t} \\
 \Rightarrow  \lim_{t \rightarrow +\infty}  \big|\tilde{\mathcal{G}}_{\alpha}( s ,t)\big| =    \lim_{t \rightarrow +\infty} e^{\frac{1}{C}\frac{k(t)}{\mu} \operatorname{Re} \{ \phi^{\alpha+1}(s) \}t}  
 \end{split}
\end{equation}
Thus, by \eqref{a}:
\begin{equation} \label{eq21*}  
\lim_{t \rightarrow +\infty}  \big|\tilde{\mathcal{G}}_{\alpha}(s ,t)\big|=0.
\end{equation}
Using the fact that $k(t)$ is bounded from bellow, \eqref{eq10*} and \eqref{a},  we have:
\begin{equation} \label{aa}  
\begin{split}
\int _{0}^{t} \tilde{ \mathcal{G}}_{\alpha}( s ,t-\tau)  d\tau = \int _{0}^{t}   e^{\frac{1}{C}\frac{k(t)}{\mu} \phi^{\alpha+1}(s)(t-\tau)}  d\tau\\
\Rightarrow  \big| \int _{0}^{t}  \tilde{\mathcal{G}}_{\alpha}(s,t-\tau)  d\tau\big|\leq \int _{0}^{t}   \big| e^{\frac{1}{C}\frac{k(t)}{\mu} \phi^{\alpha+1}(s)(t-\tau)} \big|  d\tau\\
\leq \int _{0}^{t}  e^{\frac{1}{C}\frac{|k_{0}|}{\mu} 
\operatorname{Re}  \{\phi^{\alpha+1}(s)\}(t-\tau)}d\tau\\
\Rightarrow  \lim_{t \rightarrow +\infty} \big| \int _{0}^{t}  \tilde{\mathcal{G}}_{\alpha}( s ,t-\tau)  d\tau\big| \leq  \lim_{t \rightarrow +\infty} \int _{0}^{t} e^{\frac{1}{C}\frac{|k_{0}|}{\mu} 
\operatorname{Re}  \{\phi^{\alpha+1}(s)\}(t-\tau)}d\tau \\
=\lim_{t \rightarrow +\infty}\frac{-1}{\frac{1}{C}\frac{|k_{0}|}{\mu} 
\operatorname{Re} \{\phi^{\alpha+1}(s)\}}[1-e^{\frac{1}{C}\frac{|k_{0}|}{\mu} 
\operatorname{Re}\{\phi^{\alpha+1}(s)\}t}]\\
=\frac{-1}{\frac{1}{C}\frac{|k_{0}|}{\mu} 
\operatorname{Re} \{\phi^{\alpha+1}(s)\}}
=\frac{1}{\frac{1}{C}\frac{|k_{0}|}{\mu} 
|s|^{\alpha+1} \sin(\frac{\pi}{2} \alpha)}
\end{split}
\end{equation}
 Thus, using \eqref{eq21*} and \eqref{aa}, \eqref{eq19} becomes:
\begin{equation} \label{eq2500}  
\begin{split}
\lim_{t \to +\infty} \quad \norm{ P(x, t)}_{\infty} \leq \frac{|r|}{2 \pi} \frac{1}{\frac{1}{C}\frac{|k_{0}|}{\mu} 
 \sin(\frac{\pi}{2} \alpha)} \int_{-\infty}^{+\infty} \big|\tilde{q}(s) \big| \frac{1}{|s|^{\alpha+1}}ds+\big|u_{max}\big|\\
 \end{split}
\end{equation}
Now, we focus on computing the integral in \eqref{eq2500}.  Because the term $1/|s|^{\alpha+1}$ is undefined around zero, the integral was divided to study the convergence of the integral in \eqref{eq2500} using Theorem \ref{th2.11} as follows:
\begin{equation} \label{eq2600}  
\begin{split}
 \int_{-\infty}^{+\infty} \big|\tilde{q}(s) \big| \frac{1}{|s|^{\alpha+1}}ds   = \int_{-1}^{1} \big|\tilde{q}(s) \big| \frac{1}{|s|^{\alpha+1}}ds\\
 + \int_{-\infty}^{-1} \big|\tilde{q}(s) \big| \frac{1}{|s|^{\alpha+1}}ds
  + \int_{1}^{+\infty} \big|\tilde{q}(s) \big| \frac{1}{|s|^{\alpha+1}}ds\\
 \end{split}
\end{equation}
We start by computing the following integral:
\begin{equation} \label{eq2700}  
\begin{split}
 \int_{-\infty}^{-1} \big|\tilde{q}(s) \big| \frac{1}{|s|^{\alpha+1}}ds
  \leq \norm {q(x)}_{L^{1}(\mathbb{R})} \int_{-\infty}^{-1}  \frac{1}{|s|^{\alpha+1}}ds\\
  =\norm {q(x)}_{L^{1}(\mathbb{R})} \int_{1}^{+\infty}  \frac{1}{s^{\alpha+1}}ds\\
  \Rightarrow \int_{-\infty}^{-1} \big|\tilde{q}(s) \big| \frac{1}{|s|^{\alpha+1}}ds
  \leq\frac{1}{\alpha}\norm {q(x)}_{L^{1}(\mathbb{R})}\\
 \end{split}
\end{equation}
using Theorem \ref{th2.11}, following the same reasoning as in \eqref{eq2700}, we have that:
\begin{equation} \label{eq2800}  
\begin{split}
 \int_{1}^{+\infty} \big|\tilde{q}(s) \big| \frac{1}{|s|^{\alpha+1}}ds
  \leq\frac{1}{\alpha}\norm {q(x)}_{L^{1}(\mathbb{R})}\\
 \end{split}
\end{equation}
Now, we compute the remaining integral in \eqref{eq2600}:
\begin{equation} \label{eq2900}  
\begin{split}
 \int_{-1}^{1}  \big|\tilde{q(s) }\big| \frac{1}{|s|^{\alpha+1}} ds= \int_{-1}^{1}  \big|\tilde{q(s) } \frac{1}{s^{\alpha+1}} \big| ds\\
  = \int_{-1}^{0}  \big|\tilde{q(s) } \frac{1}{s^{\alpha+1}} \big| ds+\int_{0}^{1}  \big|\tilde{q(s) } \frac{1}{s^{\alpha+1}} \big| ds\\= \lim_{\epsilon \rightarrow 0}\int_{-1}^{-\epsilon}  \big|\tilde{q(s) } \frac{1}{s^{\alpha+1}} \big| ds+\lim_{\epsilon \rightarrow 0}\int_{\epsilon}^{1}  \big|\tilde{q(s) } \frac{1}{s^{\alpha+1}}| ds\\= \lim_{\epsilon \rightarrow 0}\int_{-1}^{-\epsilon}  \big|  \mathcal{F} \{ \mathcal{J}_{y}^{\alpha+1} q(y)  \}(s) \big| ds+\lim_{\epsilon \rightarrow 0}\int_{\epsilon}^{1}  \big|  \mathcal{F} \{ \mathcal{J}_{y}^{\alpha+1} q(y)  \}(s) \big| ds\\
 \leq \lim_{\epsilon \rightarrow 0} (2-2\epsilon )\norm{   \mathcal{J}_{y}^{\alpha+1} q(y)   }_{L^{1}(\mathbb{R})}=2\norm{   \mathcal{J}_{y}^{\alpha+1} q(y)   }_{L^{1}(\mathbb{R})}\\
 \Rightarrow \int_{-1}^{1}  \big|\tilde{q(s) }\big| \frac{1}{|s|^{\alpha+1}} ds \leq 2\norm{   \mathcal{J}_{y}^{\alpha+1} q(y)   }_{L^{1}(\mathbb{R})}
 \end{split}
\end{equation}
using Theorems \ref{th2.11} and \ref{th2.12}. Then, using \eqref{eq2600}-\eqref{eq2900}, \eqref{eq2500} becomes
\begin{equation} \label{eq2300} 
\begin{split}
\lim_{t \to +\infty} \quad \norm{ P(x, t)}_{\infty} \leq  \frac{|r|}{ 2 \pi} \frac{1}{\frac{1}{C}\frac{|k_{0}|}{\mu} 
 \sin(\frac{\pi}{2} \alpha)} \big(2 \norm{   \mathcal{J}_{y}^{\alpha+1} q(y)   }_{L^{1}(\mathbb{R})}+\norm{ q(y)   }_{L^{1}(\mathbb{R})} \big)\\
\end{split}
\end{equation}
}
\end{proof}

\section {Reference Tracking process}
In this section, we study the reference tracking problem of the system given in \eqref{eq4}. We aim to track the pressure gradient at final position for the boundary controlled space FPDE in the presence of some distributions using some measurements.

Consider the space FPDE \eqref{eq4} with $L=1, g_{0}(x)=0, Q(x,t)=d_{1}(t)f(x)$, where $d_{1}(t)$ is a disturbance. We also consider a disturbance $d_{2}(t)$ on the boundary, leading to 
the following system:
\begin{equation}\label{eq21} 
\left\{\begin{array}{llll}    
 \dfrac{\partial P(x,t)}{\partial t}  - \dfrac{k(t)}{C \mu} \dfrac{\partial}{\partial x} _{}^{C}D^{\alpha}_{x}  P(x,t)=\dfrac{1}{C \rho}d_{1}(t)f(x),\\
 P(x, 0)=0,
\\
P(1, t) = u(t),   \quad  P_{x}(0, t) = d_{2}(t),   \\
\end{array}\right.
\end{equation}
where $t > 0, x \in [0,1]$, $P(x, t)$ is the gas pressure distributed in space and in time. The permeability $k(t)$ and the viscosity $\mu$ are both positive, $f(x)$ is the source term. $u(t)$ represents the control input. $_{}^{C}D^{\alpha}_{x}  P(x,t)$ is the  Caputo space fractional derivative of order $\alpha$ with $0 < \alpha \leq 1$. \R{It is important to emphasize that we are dealing with steplike and sinusoidal disturbances.} For the purpose of tracking, we consider the measurements $P(x,t)|_{x=0}$  and the output to be tracked $y(t)$ such that: $P_{x}(x,t)|_{x=1}= y(t).$
\begin{remark}
In the tracking process we will not use the measurements ($P(x,t)|_{x=0}$). Instead, we will use its fractional in space derivative of order $\alpha$. This is possible thanks to the results in \R{\cite{8483491} and \cite{8866466}}, where it has been proven that the fractional derivative of a signal can be estimated using the measurements of the signal even if this signal is noisy. We propose then to use the measurement $y_{m}(t)$ where,  $(_{}^{C}D^{\alpha}_{x}  P(x,t))|_{x=0}= y_{m}(t).$
\end{remark}
 The objective is to track asymptotically the output $y(t)$ to a desired trajectory $y_{d}(t)$, in other words we want to ensure the following : 
\begin{equation}\label{eq22} 
\lim _{t\to +\infty} e(t):=\lim _{t\to +\infty}  \big( y(t)- y_{d}(t)    \big)=0.
\end{equation}
\subsection{Output regulation using Voltera integral transformation}
In this part, we extend the well-known  backstepping approach to the space FPDE in order to design a controller that guarantees the state feedback output regulation for the considered problem. Using the same analogy as in \cite{deutscher2015backstepping}, the disturbances $d_{1}(.)$, $d_{2}(t)$, the measurements $y_{m}(t)$, and the reference $y_{d}(t)$ can be written in the space spanned by the finite dimensional signal $V(t)$ which satisfies:
\begin{equation}\label{eq23} 
\left\{\begin{array}{llll}    
V'(t)=\mathcal{S} V(t), \\
\\
V(0) = V_{0},  \\
\end{array}\right.
\end{equation}
where, $t  > 0$, $\mathcal{S}$ is a known matrix having distinct and negative eigenvalues. Thus we can write:
\begin{equation}\label{eq24} 
\left\{\begin{array}{llll}    
V'(t)=\mathcal{S} V(t),  t  > 0 \quad V(0) = V_{0},\\
\\
d_{1}(t)=a^{T}V(t),\quad d_{2}(t)=b^{T}V(t), \\
y_{d}(t)=c^{T}V(t),  y_{m}(t)=q^{T}V(t)\\
\end{array}\right.
\end{equation}

\R{where, $V_{0},a,b,c,q \in \mathbb{C}^{n_{V}}$ with $n_{V}$ an arbitrary chosen order.} \R{System \eqref{eq24} allows the modelling of unknown steplike and sinusoidal exogenous signals.} We start by introducing the Voltera coordinates transformation \cite{zhou2018boundary}, \cite{krstic2008boundary} and \cite{smyshlyaev2005backstepping}:
\begin{equation}\label{eq25} 
w(x,t)=\mathcal{V}\{P(.,t)\}(x):=P(x,t)-\int_{0}^{x}K(x,y)P(y,t)dy.
\end{equation}
This transformation is invertible, the formula for the inverse can be found in \cite{krstic2008boundary}. We consider the following space FPDE  target system:
\begin{equation}\label{eq26} 
\left\{\begin{array}{llll}    
 \dfrac{\partial w(x,t)}{\partial t}  - \dfrac{k(t)}{C \mu} \dfrac{\partial}{\partial x} _{}^{C}D^{\alpha}_{x}  w(x,t)=r^{T}(x)V(t), \\
\\
( _{}^{C}D^{\alpha}_{x}  w(x,t) )|_{x=0}=q^{T}V(t), \quad  w_{x}(0, t) =b^{T}V(t),   \\
w(1, t) = m^{T}V(t), 
\end{array}\right.
\end{equation}
where, 
\vspace{-0.5cm}
\begin{equation*} 
\begin{split}
 \quad \quad \quad  \quad r^{T}(x)=a^{T}\mathcal{V}\{f\}(x)  \\
-b^{T}\frac{k(t)}{C \mu}[_{}^{C}D^{\alpha}_{x}  K(x,y)]|_{y=0}+q^{T} \frac{k(t)}{C \mu} K(x,0),  \quad \quad  \quad \quad \quad \quad \quad \quad&
 \end{split}                   
\end{equation*}
and 
$
u(t)=\int_{0}^{1}K(1,y)P(y,t)dy+m^{T}V(t), \\
$
with the kernel system given by: 
\vspace{-0.1cm}
\begin{equation}\label{eq27} 
\left\{\begin{array}{llll}    
\frac{\partial}{\partial x}\R{\bigg(} \quad \hspace{-4mm} _{y}^{C}D^{\alpha}_{x} K(x,y)\R{\bigg)}=_{}^{C}D^{\alpha}_{y,x} \R{\bigg(} \frac{\partial}{\partial y} K(x,y)\R{\bigg)}, \\
_{y}^{C}D^{\alpha-1}_{x} K(x,y)|_{y=x}=0,\\
_{}^{C}D^{\alpha}_{y,x} K(x,y)|_{y=x}=0,\\
k(x,x)=0, \\
 \R{K(x,0) \neq 0}
\end{array}\right.
\end{equation}
where, ${}^{C}D^{\alpha}_{y,x}$ is the right sided Caputo derivative defined in \eqref{sided}

\begin{lemma}
Using transformation \eqref{eq25}, if there exists a twice continuously differentiable kernel function $K(x, y)$ satisfying \eqref{eq27} then, system \eqref{eq21} is equivalent to \eqref{eq26}.
\end{lemma}
\vspace{-0.5cm}
\begin{proof}
We start by computing the time classical derivative of the new coordinate $w(x, t)$:
\vspace{-0.4cm}
\begin{equation}\label{eq28} 
\begin{split}
\frac{\partial}{\partial t} w(x, t)=  \frac{\partial}{\partial t} P(x,t)-\int_{0}^{x}k(x,y) \frac{\partial}{\partial t} P(y,t)dy \\
= \frac{k(t)}{C \mu} \frac{\partial}{\partial x} _{}^{C}D^{\alpha}_{x} P(x,t)-\frac {k(t)}{C \mu} \int_{0}^{x}K(x,y)\frac{\partial}{\partial y}  ^{C}D^{\alpha}_{y}P(y,t) dy\\
 +\frac{1}{C \rho}a^{T}V(t)f(x) -\frac{1}{C \rho}\int_{0}^{x}k(x,y)a^{T}V(t)f(y)dy \\
=\frac {k(t)}{C \mu} \frac{\partial}{\partial x}  _{}^{C}D^{\alpha}_{x}P(x,t)+\frac{1}{C \rho}a^{T}V(t)\mathcal{V} \{f\}(x)\\
 +  \frac{k(t)}{C \mu}  \int_{0}^{x} \frac{\partial}{\partial y} K(x,y) _{}^{C}D^{\alpha}_{y} P(y,t) dy \\
  - \frac{k(t)}{C \mu} [K(x,y) _{}^{C}D^{\alpha}_{y} P(y,t)]_{y=0}^{y=x}\\
=\frac {k(t)}{C \mu} \frac{\partial}{\partial x} _{}^{C}D^{\alpha}_{x}P(x,t)+\frac{1}{C \rho}a^{T}V(t)f(x) \\
 - \frac{k(t)}{C \mu}  \int_{0}^{x} \hspace{1mm}  _{}^{C}D^{\alpha}_{y,x}  \frac{\partial}{\partial y} K(x,y) P(y,t) dy  \quad \\
+ \frac{k(t)}{C \mu}  [    _{}^{C}D^{\alpha-1}_{y,x}   \frac{\partial}{\partial y} K(x,y)  P(y,t) ]_{y=0}^{y=x}\\
- \frac{k(t)}{C \mu} [K(x,y) _{}^{C}D^{\alpha}_{y}  P(y,t)]_{y=0}^{y=x}\\
 -\frac{1}{C \rho} a^{T}V(t) \int_{0}^{x}K(x,y)f(y)dy 
 \end{split}                                      
\end{equation}

where, \eqref{eq28} is obtained by first applying a classical integration by parts then the fractional integration by parts \eqref{eq2}. 
Using the generalized Leibnitz differentiation rule \eqref{eq3}, we obtain the Caputo spatial fractional derivative of the new coordinate $w(x, t)$:
\begin{equation}\label{eq29} 
\begin{split}
_{}^{C}D^{\alpha}_{x}  w(x, t)= _{}^{C}D^{\alpha}_{x} P(x,t) -\int_{0}^{x} \hspace{1mm} _{y}^{C}D^{\alpha}_{x} K(x,y) P(y,t)dy \\
-[_{y}^{C}D^{\alpha-1}_{x} K(x,y)P(y, t)]|_{y=x}, \\
\end{split}
\end{equation}
then, by the classical Leibnitz differentiation rule, we obtain:

\begin{equation}\label{eq30} 
\begin{split}
\frac{\partial}{\partial x}   \hspace{1mm} _{}^{C}D^{\alpha}_{x} w(x, t)= \frac{\partial}{\partial x}  _{}^{C}D^{\alpha}_{x} P(x,t)  \\
-\int_{0}^{x} \frac{\partial}{\partial x}   _{y}^{C}D^{\alpha}_{x} K(x,y) P(y,t)dy\\
- [_{y}^{C}D^{\alpha}_{x} K(x,y)P(y, t)]|_{y=x}\\
-\frac{\partial}{\partial x} ( [_{y}^{C}D^{\alpha-1}_{x} K(x,y)P(y, t)]|_{y=x}).\\
\end{split}
\end{equation}

Thus, by \eqref{eq21}, \eqref{eq24}, \eqref{eq25}, \eqref{eq27}, \eqref{eq28}, \eqref{eq29} and \eqref{eq30} we have:

\begin{equation}\label{eq31} 
\begin{split}
\frac{\partial w(x,t)}{\partial t} \frac{k(t)}{C \mu} \frac{\partial}{\partial x} _{}^{C}D^{\alpha}_{x} w(x,t)\\
= \quad \bigg(\frac{1}{C \rho} a^{T}\mathcal{V}\{f\}(x)-b^{T}\frac{k(t)}{C \mu}  [_{}^{C}D^{\alpha-1}_{y,x}  \frac{\partial}{\partial y} K(x,y)  ]_{y=0} \\
 +q^{T} \frac{k(t)}{C \mu} K(x,0)\bigg) V(t) \\
\end{split}  
\end{equation}

By looking to \eqref{eq31}, we conclude that we need an extra condition on the kernel and which is $ K(x,0) \neq 0$ in order not to loss the measurements. 
We set $ w(1, t)=w_{1}(t)$., where by \eqref{eq24}, $w_{1}(t)$ can be written as: $w_{1}(t)=m^{T}V(t)$

Using \eqref{eq21} and transformation \eqref{eq25}, we obtain :
\begin{equation}\label{eq32} 
u(t)=P(1, t)=\int_{0}^{1}K(1,y)P(y,t)dy+m^{T}V(t).
\end{equation}
Which completes the proof.
\end{proof}

\begin{lemma}
The kernel system \eqref{eq27} admits at least one family of twice continuously differentiable solutions $K(x,y)$ in the triangle $ 0 \leq y \leq x \leq 1$ and which is given by:.
\begin{equation}\label{eq33} 
K(x,y)=(x-y)^{2m+1}, \quad m \in  \mathbb{N} 
\end{equation}
\end{lemma}

\begin{proof}
Consider system \eqref{eq27} in $ 0 \leq y \leq x \leq 1$. 
To check that the kernel function given in \eqref{eq33} is valid, we start by computing both sides of the kernel PDE using \R{T}heorem \ref{thm2}:
\begin{equation}\label{eq34} 
\frac{\partial}{\partial x} _{y}^{C}D^{\alpha}_{x} (x-y)^{2m+1}=\frac{\Gamma(2m+2)}{\Gamma(2m+1-\alpha)}(x-y)^{2m-\alpha},
\end{equation}
and
\begin{equation}\label{eq35} 
_{}^{C}D^{\alpha}_{y,x} \frac{\partial}{\partial y} (x-y)^{2m+1}=\frac{\Gamma(2m+2)}{\Gamma(2m+1-\alpha)}(x-y)^{2m-\alpha},
\end{equation}
By \eqref{eq34} and \eqref{eq35} the proposed kernel function satisfies the kernel PDE. Let's check the boundary conditions:
we have:
$$
_{y}^{C}D^{\alpha-1}_{x} K(x,y)|_{y=x}=\frac{\Gamma(2m+2)}{\Gamma(2m+3-\alpha)}(x-y)^{2m+2-\alpha}|_{y=x}=0,\\
$$
which means that:
$$
_{}^{C}D^{\alpha}_{y,x} K(x,y)|_{y=x}=\frac{\Gamma(2m+2)}{\Gamma(2m+2-\alpha)}(x-y)^{2m+1-\alpha}|_{y=x}=0,\\
 $$
 and 
$$
K(x,x)=0.
$$
\R{Finally, $ K(x,0) =(x)^{2m+1}\neq 0$,} which completes the proof.
\end{proof}

The objective now is to determine $m^{T}$ which guarantees that $w(x,t)$ will achieve the output regulation at steady state. We define the tracking error:
\begin{equation}\label{eq36}  
e(x,t)=w(x,t)-M^{T}(x)V(t), \\
\end{equation}
where $M^{T}$ has to be determined. Let's now define the following systems:
\begin{equation}\label{eq37} 
\left\{\begin{array}{llll}    
\frac{\partial e(x,t)}{\partial t}  = \frac{k(t)}{C \mu} \frac{\partial}{\partial x} _{}^{C}D^{\alpha}_{x} e(x,t) \\
\\
 e_{x}(0, t) =0,  \quad e(1, t) =0,   \\
 \\
\big( \mathcal{V}^{-1}\{ e_{x}(x, t) \} \big)|_{x=1} =e(t), \\
\end{array}\right.
\end{equation}
and the system:
\begin{equation}\label{eq38} 
\left\{\begin{array}{llll}    
M^{T}(x)  \mathcal{S}-\frac{k(t)}{C \mu} \frac{\partial}{\partial x} _{}^{C}D^{\alpha}_{x} M^{T}(x)=r^{T}(x), \\
\\
M_{x}^{T}(x)|_{x=0}=b^{T},   \\
\big( \mathcal{V}^{-1}\{ M^{T}(x) \} \big)|_{x=1} =c^{T}, \\
\end{array}\right.
\end{equation}
\\
\begin{theorem}
Using \eqref{eq23}, \eqref{eq25} and \eqref{eq36} if there exists a  \R{twice continuously differentiable function $M^{T}(x)$ } solution of \eqref{eq38} then, the tracking error $e(x,t)$ satisfies \eqref{eq37}\R{,}
where, $m^{T}$ is chosen \R{as follows:}
$$
[M^{T}(x)]|_{x=1}=m^{T},
$$
\end{theorem}
 \begin{proof}
 We start by computing the time classical derivative and the space fractional derivative of the tracking error \eqref{eq36}, and replace them in \eqref{eq21} we get:
 \begin{equation}\label{eq39} 
\begin{split}  
\frac{\partial e(x,t)}{\partial t}-\frac{k(t)}{C \mu} \frac{\partial}{\partial x} _{}^{C}D^{\alpha}_{x}e(x,t)= \qquad  \qquad  \qquad  \\ 
 r^{T}(x)V(t)-M^{T}(x)  \mathcal{S}V(t)+\frac{k(t)}{C \mu} \frac{\partial}{\partial x} _{}^{C}D^{\alpha}_{x} M^{T}(x)V(t). \\
\end{split}
\end{equation} 
We take $M^{T}(x)$ to be solution of \eqref{eq38}. Thus, \eqref{eq39} becomes \eqref{eq37}. Furthermore if we chose $m^{T}$ that satisfies:
$$
[M^{T}(x)]|_{x=1}=m^{T},
$$
which is in \eqref{eq37} equivalent to the condition:
$$
e(1, t) =0.
$$
The control in \eqref{eq32} becomes:
\begin{equation}\label{eq40} 
u(t)=P(1, t)=\int_{0}^{1}K(1,y)P(y,t)dy+[M^{T}(x)]|_{x=1}V(t).
\end{equation}
 \end{proof}

\begin{theorem}\label{error}
The tracking error system \eqref{eq37} is asymptotically stable \R{in $\B{L_{\infty}(\mathbb{R})}$}.
\end{theorem}

\begin{proof}
The tracking error function given by \eqref{eq37} is a particular case of \eqref{eq4} for a source term equal to zero. Thus, the solution of \eqref{eq37} is given in \R{T}heorem \ref{thm1} by \eqref{eq13}:
\begin{equation} \label{eq41}  
\begin{split}
e(x, t)= \int _{- \infty}^{+ \infty}  \mathcal{G}_{\alpha}(\mid x-y \mid,t)e_{0}(y) dy.  \\
=\frac{1}{2\pi}\int _{- \infty}^{+ \infty} e^{-isx}  \tilde{\mathcal{G}}_{\alpha}(s,t)e_{0}(s) ds.
\end{split}
\end{equation}
where,
\begin{equation} \label{b} 
\begin{split}
e_{0}(x)=e(x, 0)=g_{0}(x)-&\int_{0}^{x}K(x,y)g_{0}(y)dy-M^{T}(x)V_{0}\\
=-M^{T}(x)V_{0},
\end{split}
\end{equation} 
\R{
 $M^{T}(x)$ is twice  continuously differentiable, thus, $e_{0}(x)$ is in $L^{1}(\mathbb{R})$ which means that its Fourier transform  $\tilde{e}_{0}(s)$ exist. Then, we get:
\begin{equation} \label{eq17b}  
\begin{split}
\big|e(x, t)\big|   
\leq \frac{1}{2 \pi} \int _{-\infty}^{+\infty}  \big|   \tilde{\mathcal{G}}_{\alpha}(s,t) e^{-isx}\tilde{e}_{0}(s)  \big|  ds\\
= \frac{1}{2\pi}\int _{-\infty}^{+\infty}  \big| \tilde{\mathcal{G}}_{\alpha}(s,t)\big| \big|\tilde{e}_{0}(s) \big|ds.  \\
\end{split}
\end{equation}
By taking the sup on $x$ over $[0,L]$ of \eqref{eq17b}, we get:
\begin{equation} \label{eq18b}  
\begin{split}
 \norm{ e(x, t) }_{\infty} \leq 
\frac{1}{2 \pi}  \int _{-\infty}^{+\infty}  \big|\tilde{\mathcal{G}}_{\alpha}(s,t)\big| \big| \tilde{e}_{0}(s) \big| ds,  \\
\end{split}
\end{equation}
taking the limit of \eqref{eq18b}, we get:
\begin{equation} \label{eq19b}  
\begin{split}
\lim_{t \rightarrow +\infty} \norm{ e(x, t) }_{\infty} \leq  \frac{1}{2 \pi} \int _{-\infty}^{+\infty}   \lim_{t \rightarrow +\infty}   \big|\tilde{\mathcal{G}}_{\alpha}(s,t)\big| \big| \tilde{e}_{0}(s) \big| ds,  \\
\end{split}
\end{equation}
where
\begin{equation*}   
\tilde{\mathcal{G}}_{\alpha}(s,t)= e^{\frac{1}{C}\frac{k(t)}{\mu} \phi^{\alpha+1}(s)t}
\end{equation*}
and
$$\phi^{\alpha+1}(s)=(-is)^{\alpha+1}.  $$\\
 On the other hand, we have from equation \eqref{eq21*}  that:
 \begin{equation} \label{eq21*b}  
\lim_{t \rightarrow +\infty}  \big|\tilde{\mathcal{G}}_{\alpha}(s ,t)\big|=0,
\end{equation}
using the fact that both the permeability and the viscosity are both positive and the linear-compressibility with $C > 0$. Thus, using \eqref{eq21*b}, \eqref{eq19b} becomes:
\begin{equation} \label{eq20b}  
\begin{split}
\lim_{t \rightarrow +\infty} \norm{ e(x, t) }_{\infty} =0  \\
\end{split}
\end{equation}
which completes the proof.}
 \end{proof}
 
 \begin{remark}
 Notice that the kernel system in \eqref{eq27} is non-integer order PDE with very complex boundary conditions and with right sided fractional derivative $_{}^{C}D^{\alpha}_{y,x}$ and shifted fractional derivative $_{y}^{C}D^{\alpha}_{x}$ which makes it hard to solve. In \eqref{eq33}, we propose one possible family of solutions which satisfies \eqref{eq27}. To avoid having this non-integer kernel PDE, in the next part, we propose an alternative method which simplifies the Kernel PDE into an integer PDE which has a unique solution.

  \end{remark}
 
 %%%%%%%%%%%%%%%%%%%%%%%%%%%%%%%%%%%%%%%%%%%%%%%%%%%%%%%%%%%%%%%%%%%%%%%%%%%%%%%%%%%%

\subsection{Novel coordinates transformation for output regulation}
In this part, we propose an alternative coordinates transformation to the Voltera transformation \eqref{eq25} which allows to achieve a stable target system similar to \eqref{eq26}  but with an integer kernel system simpler than the one in \eqref{eq27} and with less complex derivations. This method can be used for PDEs with a more complex space operator than the time one (which is the case for our system \eqref{eq21}).

We assume that the disturbances $d_{1}(.)$, $d_{2}(t)$, the measurements $y_{m}(t)$, and the reference $y_{d}(t)$ in \eqref{eq21} can be written in the space spanned by the finite dimensional signal $V(t)$ which satisfies \eqref{eq23} and \eqref{eq24}.
We start by introducing the new coordinates transformation :
\begin{equation}\label{eq25N} 
w(x,t)=\mathcal{V}\{P(.,t)\}(x):=P(x,t)-\int_{0}^{t}P(x,\tau)l(\tau,t)d\tau.
\end{equation}
This transformation is invertible, the formula for the inverse will be discussed later. We consider the following space fractional target system:
\begin{equation}\label{eq26N} 
\left\{\begin{array}{llll}    
 \dfrac{\partial w(x,t)}{\partial t}  - \frac{k(t)}{C \mu} \dfrac{\partial}{\partial x} _{}^{C}D^{\alpha}_{x}  w(x,t)=\dfrac{1}{C \rho}a^{T}\mathcal{V}\{f\}(x)V(t), \\
\\
w(0, t) =n^{T}V(t) ,   \\
w(1, t) = m^{T} V(t),   \\
\end{array}\right.
\end{equation}
where, 
$$
u(t)=\int_{0}^{t}P(1,\tau)l(\tau,t)d\tau+m^{T}V(t).
$$
and the kernel system:
\begin{equation}\label{eq27N} 
\left\{\begin{array}{llll}    
\dfrac{\partial}{\partial t}  l(\tau,t)=- \dfrac{\partial}{\partial \tau} l(\tau,t), \\
l(0,t)=\hat{l}(t),\\
\end{array}\right.
\end{equation}
where, $\hat{l}(.)$ is a known differentiable function. 
\begin{theorem}\label{replace}
Using transformation \eqref{eq25N}, if there exist a continuously differentiable kernel function $l(t, \tau)$ satisfying \eqref{eq27N}  then, system \eqref{eq21} is equivalent to \eqref{eq26N}.
\end{theorem}

\begin{proof}
We start by computing the time classical derivative of the new coordinate $w(x, t)$:
\begin{equation}\label{eq28N} 
\begin{split}
\frac{\partial}{\partial t} w(x, t)=  \frac{\partial}{\partial t} P(x,t) \qquad \qquad\\
 -\int_{0}^{t}P(x,\tau) \frac{\partial}{\partial t} l(\tau,t)d \tau -P(x,t) l(t,t),
\end{split}
\end{equation}
 using the non-integer Leibnitz rule. The Caputo spatial fractional derivative of the new coordinate $w(x, t)$ is given by:
 \vspace{-3mm}
\begin{equation}\label{eq29N} 
_{}^{C}D^{\alpha}_{x} w(x, t)= \quad \hspace{-4mm} _{}^{C}D^{\alpha}_{x}P(x,t)-\int_{0}^{t} \hspace{1mm} _{}^{C}D^{\alpha}_{x} P(x,\tau)l(\tau,t)d \tau
\end{equation}
Thus, we have:
\begin{equation}\label{eq30N} 
\begin{split}
\frac{k(t)}{C \mu} \frac{\partial}{\partial x} \quad \hspace{-4mm}_{}^{C}D^{\alpha}_{x} w(x, t)=  \frac{k(t)}{C \mu} \frac{\partial}{\partial x} \quad \hspace{-4mm} _{}^{C}D^{\alpha}_{x} P(x,t)\\
-\int_{0}^{t} \frac{k(t)}{C \mu} \frac{\partial}{\partial x}\quad \hspace{-4mm}_{}^{C}D^{\alpha}_{x} P(x,\tau)l(\tau,t) d \tau \\
=  \frac{\partial}{\partial t} P(x,t)-
 \int_{0}^{t} \frac{\partial}{\partial \tau} P(x,\tau) l(\tau,t) d \tau\\
 -\frac{1}{C \rho}a^{T}V(t)f(x) +\frac{1}{C \rho}\int_{0}^{t}l(\tau,t)a^{T}V(\tau)f(x)d\tau \\
 =\frac{\partial}{\partial t} P(x,t)+
 \int_{0}^{t}  P(x,\tau) \frac{\partial}{\partial \tau} l(\tau,t) d \tau -[ P(x,\tau) l(\tau,t)]|_{\tau=0}^{\tau=t}\\
 -\frac{1}{C \rho}a^{T}V(t)f(x) +\frac{1}{C \rho}\int_{0}^{t}l(\tau,t)a^{T}V(\tau)f(x)d\tau \\
 \end{split}                                      
\end{equation}
where, \eqref{eq30N} is obtained by applying a classical integration by parts. Thus, by \eqref{eq21}, \eqref{eq24}, \eqref{eq25N}, \eqref{eq28N}, and \eqref{eq30N} we have:
\vspace{-3mm}
\begin{equation}\label{eq31N} 
\begin{split}
\frac{\partial w(x,t)}{\partial t} \frac{k(t)}{C \mu} \frac{\partial}{\partial x} \quad \hspace{-4mm}_{}^{C}D^{\alpha}_{x} w(x,t)\\
= \quad \frac{1}{C \rho} a^{T}\mathcal{V}\{f\}(x) V(t) \\
\end{split}  
\end{equation}

We set $ w(1, t)=w_{1}(t)$, where by \eqref{eq24}, $w_{1}(t)$ can be written as: $w_{1}(t)=m^{T}V(t)$

Using \eqref{eq21} and transformation \eqref{eq25N}, we obtain :
\begin{equation}\label{eq32N} 
u(t)=P(1, t)=\int_{0}^{t}P(1,\tau)l(\tau,t)d\tau+m^{T}V(t).
\end{equation}
Which completes the proof.
\end{proof}
\begin{theorem}
The kernel system \eqref{eq27N} admits a unique continuously differentiable solution $l(t,\tau)$ in the triangle $ 0 \leq \tau \leq t $ and which is given by:
\vspace{-3mm}
\begin{equation}\label{eq33N} 
l(\tau,t)=\hat{l}(t-\tau), 
\end{equation}
\vspace{-6mm}
\end{theorem}
The objective now is to determine $m^{T}$ which guarantees that $w(x,t)$ will achieve the output regulation at steady state. We define the tracking error:
\begin{equation}\label{eq36N}  
E(x,t)=w(x,t)-M^{T}(x)V(t), \\
\end{equation}
where $M^{T}$ has to be determined. Let's now define the following systems:
\begin{equation}\label{eq37N} 
\left\{\begin{array}{llll}    
\dfrac{\partial E(x,t)}{\partial t}  = \frac{k(t)}{C \mu} \dfrac{\partial}{\partial x} \quad \hspace{-4mm} _{}^{C}D^{\alpha}_{x} E(x,t) \\
\\
 E(0, t) =0,  \quad E(1, t) =0,   \\
 \\
\big( \mathcal{V}^{-1}\{ E_{x}(x, t) \} \big)|_{x=1} =E(t), \\
\end{array}\right.
\end{equation}
and the system:
\begin{equation}\label{eq38N} 
\left\{\begin{array}{llll}    
M^{T}(x)  \mathcal{S}-\frac{k(t)}{C \mu} \frac{\partial}{\partial x}\quad \hspace{-4mm} _{}^{C}D^{\alpha}_{x} M^{T}(x)=r^{T}(x), \\
\\
M^{T}(x)|_{x=0}=n^{T},   \\
\big( \mathcal{V}^{-1}\{ M^{T}(x) \} \big)|_{x=1} =c^{T}, \\
\end{array}\right.
\end{equation}
\\
\begin{theorem}
Using \eqref{eq23}, \eqref{eq25N} and \eqref{eq36N} if there exists a solution of \eqref{eq38N} then, the tracking error $e(x,t)$ satisfies \eqref{eq37N}\R{,}
where, $m^{T}$ is chosen \R{as follows:}
$$
[M^{T}(x)]|_{x=1}=m^{T},
$$
\end{theorem}
 \begin{proof}
 We start by computing the time classical derivative and the space fractional derivative of the tracking error \eqref{eq36N}, and replace them in \eqref{eq21} we get:
 \begin{equation}\label{eq39N} 
\begin{split}  
\frac{\partial E(x,t)}{\partial t}-\frac{k(t)}{C \mu} \frac{\partial}{\partial x}\quad \hspace{-4mm} _{}^{C}D^{\alpha}_{x}E(x,t)= \qquad \qquad \quad  \\
r^{T}(x)V(t)-M^{T}(x)  \mathcal{S}V(t)+\frac{k(t)}{C \mu} \frac{\partial}{\partial x}\quad \hspace{-4mm} _{}^{C}D^{\alpha}_{x} M^{T}(x)V(t). \\
\end{split}
\end{equation} 
We take $M^{T}(x)$ to be solution of \eqref{eq38N}. Thus, \eqref{eq39N} becomes \eqref{eq37N}. Furthermore if we chose $m^{T}$ that satisfies:
$$
[M^{T}(x)]|_{x=1}=m^{T},
$$
which is in \eqref{eq37N} equivalent to the condition:
$$
E(1, t) =0.
$$
The control in \eqref{eq32N} becomes:
\begin{equation}\label{eq40N} 
u(t)=P(1, t)=\int_{0}^{t}P(1,\tau)l(\tau,t)d\tau+[M^{T}(x)]|_{x=1}V(t).
\end{equation}
 \end{proof}

\begin{theorem}
The tracking error system \eqref{eq37N} is asymptotically stable \R{in $\B{L_{\infty}(\mathbb{R})}$}.
\end{theorem}
The proof is similar to the proof of \R{T}heorem \ref{error}
\begin{theorem}
The inverse of transformation \eqref{eq25N} is given by:
\begin{equation}\label{eq41N} 
P(x,t)=\mathcal{V}^{-1}\{w(.,t)\}(x):=w(x,t)+\int_{0}^{t}w(x,\tau)L(\tau,t)d\tau.
\end{equation} 
 with:
 \begin{equation}\label{eq42N} 
 L(\tau,t)=\hat{L}(t-\tau), 
\end{equation} 
 where, 
 $$
 \hat{L}(t)=L(0,t).
 $$
\end{theorem}
 \begin{proof}
 The proof is similar to the proof of \R{T}heorem \ref{replace}. using transformation \eqref{eq42N}, we compute the first temporal derivative and the fractional spatial derivative of $P(x,t)$ and inject it in \eqref{eq21}. Using \eqref{eq21}, we conclude the following system for the kernel $L(\tau,t)$:
 \begin{equation}\label{eq43N} 
\left\{\begin{array}{llll}    
\frac{\partial}{\partial t}  L(\tau,t)=- \frac{\partial}{\partial \tau} L(\tau,t), \\
L(0,t)=\hat{L}(t),\\
\end{array}\right.
\end{equation}
\R{System \eqref{eq43N} admits a unique solution given in \cite{evans1998partial} by:
 \begin{equation*}
 L(\tau,t)=\hat{L}(t-\tau), 
\end{equation*} 
 where, 
 $$
 \hat{L}(t)=L(0,t).
 $$
 which completes the proof.
}
 \end{proof}

 \section {Adaptive boundary observer for space fractional PDEs subject to domain and boundary disturbances\\}
\R{
To design the tracking controllers in \eqref{eq40} and \eqref{eq40N}, we need to recover the state $P(x,t)$ of system \eqref{eq21}.  We recall that, it has been mentioned in Remark 6.1 that 
in the tracking process, we did not use the measurements ($z_{m}(t)=P(x,t)|_{x=0}$). Instead, we  used its fractional in space derivative of order $\alpha$ ($ y_{m}(t)=_{}^{C}D^{\alpha}_{x}  P(x,t)|_{x=0}.$). This was possible thanks to the results in [41, 42], where it has been proven that the fractional derivative of a signal can be estimated using the measurements of the signal even if this signal is noisy.
In the adaptive observer design, we propose then to use both the measurements $z_{m}(t)$ and  $y_{m}(t)$ where, $  P(x,t)|_{x=0}= z_{m}(t)$ and $(_{}^{C}D^{\alpha}_{x}  P(x,t))|_{x=0}= y_{m}(t).$ which will prove to be judicious later. Thus, we design the following observer for system \eqref{eq21}:}
\R{
 \begin{equation}\label{eq1A}
\left\{\begin{array}{llll}    
 \frac{\partial \hat{P}(x,t)}{\partial t}  - \frac{k(t)}{C \mu} \frac{\partial}{\partial x} _{}^{C}D^{\alpha}_{x}  \hat{P}(x,t)=\frac{1}{C \rho}\hat{d}_{1}(t)f(x)+u(x,t)\\
 \qquad \qquad -\frac{k(t)}{C \mu}H_{1}(x)(\hat{P}(0,t)-z_{m}(t))-\frac{k(t)}{C \mu}H_{2}(x)(_{}^{C}D^{\alpha}_{x} \hat{P}(0,t)-y_{m}(t)),\\
 \hat{P}(x, 0)=\hat{P}_{0}(x), \hat{P}(1, t) = u(t),  \\
  \hat{P}_{x}(0, t) = \hat{d}_{2}(t).   \\
\end{array}\right.
\end{equation}
} \R{
 This observer is designed similarly as in [15-17] for a parabolic systems. Where, $H_{1}(x)$ and $H_{2}(x)$ are a space-dependent observer gains, $\hat{P}_{0}(x)$ is an arbitrary initial condition  which satisfies some conditions that will be determined later.  $\hat{d}_{1}(t)$ and $\hat{d}_{2}(t)$ are the disturbances estimates and $u(x,t)$ is an additional feedback term that will be determined latter. We introduce the following usual estimation errors:\\} \R{
 \begin{equation}\label{stateer}
  \text{State estimation error: \quad}   \tilde{P}(x, t)=\hat{P}(x, t)-P(x, t). 
 \end{equation}
  \begin{equation}\label{diser}
  \text{Disturbances estimation error: \quad}   \tilde{\theta}(t)=\hat{\theta}(t)-\theta(t):=\begin{bmatrix}
\tilde{d}_{1}(t) \\
\tilde{d}_{2}(t)
\end{bmatrix},
 \end{equation}
 where, 
   \begin{equation}\label{diserr}
  \hat{\theta}(t) :=\begin{bmatrix}
\hat{d}_{1}(t) \\
\hat{d}_{2}(t)
\end{bmatrix}
\text{\quad and \quad }
\begin{bmatrix}
\tilde{d}_{1}(t) \\
\tilde{d}_{2}(t)
\end{bmatrix}:=\begin{bmatrix}
\hat{d}_{1}(t)-d_{1}(t) \\
\hat{d}_{2}(t)-d_{2}(t)
\end{bmatrix}.
 \end{equation}}
 \R{
Then, using \eqref{eq1A}, it follows that the state estimation error in \eqref{stateer} satisfies the following system: 
  \begin{equation}\label{eq2A}
\left\{\begin{array}{llll}    
 \frac{\partial \tilde{P}(x,t)}{\partial t}  - \frac{k(t)}{C \mu} \frac{\partial}{\partial x} _{}^{C}D^{\alpha}_{x}  \tilde{P}(x,t)=\frac{1}{C \rho}\tilde{d}_{1}(t)f(x)+u(x,t)\\
 \qquad \qquad -\frac{k(t)}{C \mu}H_{1}(x)\tilde{P}(0,t)-\frac{k(t)}{C \mu}H_{2}(x)_{}^{C}D^{\alpha}_{x}\tilde{P}(0,t) ,\\
 \tilde{P}(x, 0)=\hat{P}_{0}(x), \tilde{P}(1, t) =0,  \\
  \tilde{P}_{x}(0, t) =\tilde{d}_{2}(t),   \\
\end{array}\right.
\end{equation}
}
\subsection{Finite-dimensional backstepping-like transformation}\leavevmode\par\noindent
\R{The objective now is to cancel the terms on $\tilde{d}_{1}(t)$ and $\tilde{d}_{2}(t)$ in the observer error system \eqref{eq2A}. To this end, we consider the finite-dimensional backstepping-like transformation defined in [15-17] as:
\begin{equation}\label{like} 
w(x,t)=\tilde{P}(x,t)-\lambda_{1}(x,t)\tilde{d}_{1}(t)-\lambda_{2}(x,t)\tilde{d}_{2}(t):=\tilde{P}(x,t)-\Lambda(x,t)\tilde{\theta}(t),
\end{equation}
where, 
\begin{equation*}
\Lambda(x,t)=\begin{bmatrix}
\lambda_{1}(x,t) \quad \lambda_{2}(x,t).
\end{bmatrix}
\end{equation*}
Using the PDE in \eqref{eq2A}, $w(x,t)$ satisfies the following PDE:
\begin{equation}\label{eq2AA} 
\begin{split}
 \frac{\partial w(x,t)}{\partial t}  & =\frac{k(t)}{C \mu} \frac{\partial}{\partial x} _{}^{C}D^{\alpha}_{x} \tilde{P}(x,t)+\frac{1}{C \rho}\tilde{d}_{1}(t)f(x)-H_{1}(x)\tilde{P}_{x}(1,t)\\
 & +u(x,t)-\frac{\partial }{\partial t}(\Lambda(x,t)\tilde{\theta}(t)),\\
 &=\frac{k(t)}{C \mu} \frac{\partial}{\partial x} _{}^{C}D^{\alpha}_{x} w(x,t)-H_{1}(x)w_{x}(1,t)+u(x,t)-\Lambda(x,t)\dot{\tilde{\theta}}(t)\\
 &+\bigg[\frac{k(t)}{C \mu} \frac{\partial}{\partial x} _{}^{C}D^{\alpha}_{x}\lambda_{1}(x,t)-\frac{\partial }{\partial t}\lambda_{1}(x,t)-H_{1}(x)\frac{\partial }{\partial x}\lambda_{1}(x,t)+\frac{1}{C \rho}f(x) \bigg]\tilde{d}_{1}(t)\\
 &+\bigg[\frac{k(t)}{C \mu} \frac{\partial}{\partial x} _{}^{C}D^{\alpha}_{x}\lambda_{2}(x,t)-\frac{\partial }{\partial t}\lambda_{2}(x,t)-H_{1}(x)\frac{\partial }{\partial x}\lambda_{2}(x,t) \bigg]\tilde{d}_{2}(t),\\
 \end{split}
\end{equation}
This suggests the following choice of the feedback expression for $u(x,t)$:
$$u(x,t)=\Lambda(x,t)\dot{\tilde{\theta}}(t)$$ and the following trajectory of the auxiliary states:
\begin{equation}\label{FF} 
\begin{split}
\frac{k(t)}{C \mu} \frac{\partial}{\partial x} _{}^{C}D^{\alpha}_{x}\lambda_{1}(x,t)-\frac{\partial }{\partial t}\lambda_{1}(x,t)-H_{1}(x)\frac{\partial }{\partial x}\lambda_{1}(x,t)+\frac{1}{C \rho}f(x)=0\\
\frac{k(t)}{C \mu} \frac{\partial}{\partial x} _{}^{C}D^{\alpha}_{x}\lambda_{2}(x,t)-\frac{\partial }{\partial t}\lambda_{2}(x,t)-H_{1}(x)\frac{\partial }{\partial x}\lambda_{2}(x,t)=0
\end{split}
\end{equation}
Thus, \eqref{eq2AA} becomes:
\begin{equation}\label{eq2AAA} 
 \frac{\partial w(x,t)}{\partial t} 
 =\frac{k(t)}{C \mu} \frac{\partial}{\partial x} _{}^{C}D^{\alpha}_{x} w(x,t)-H_{1}(x)w_{x}(1,t)\\
\end{equation}
This is completed by the following initial and boundary conditions :
 \begin{equation}\label{F1}
\left\{\begin{array}{llll}    
\lambda_{1}(x,0)=0\\
\frac{\partial }{\partial x}\lambda_{1}(x,t)|_{x=0}=\lambda_{1}(1,t)=0\\
\end{array}\right. \quad
\left\{\begin{array}{llll}    
\lambda_{2}(x,0)=0\\
\lambda_{2}(1,t)=0\\
\frac{\partial }{\partial x}\lambda_{2}(x,t)|_{x=0}=1.\\
\end{array}\right.
\end{equation}
\begin{theorem}\label{00000}
There exist a unique auxiliary states $\lambda_{1}(x,t)$ and $\lambda_{2}(x,t)$ that satisfies the PDEs in \eqref{FF} and the initial and boundary conditions in \eqref{F1}. Furthermore, Theses auxiliary states are stable. That is:
\begin{equation*} 
\lim_{t \rightarrow +\infty} \norm{\Lambda(x, t)}_{\infty}< \infty 
\end{equation*}
\end{theorem}
\begin{proof}
The proof is given in appendix B.
\end{proof}
The choice of  initial and boundary conditions in \eqref{F1} is efficient because they allow the following initial and boundary condition for the system $w(x,t)$:
 \begin{equation}\label{F2}
\left\{\begin{array}{llll}    
w(x,0)=\hat{P}_{0}(x)\\
w_{x}(0,t)=w(1,t)=0\\
\end{array}\right.
\end{equation}
Thus the system on $w(x,t)$ become:
 \begin{equation}\label{w}
\left\{\begin{array}{llll}  
 \frac{\partial w(x,t)}{\partial t} 
 =\frac{k(t)}{C \mu} \frac{\partial}{\partial x} _{}^{C}D^{\alpha}_{x} w(x,t)-H_{1}(x)w(0,t)-H_{2}(x)_{}^{C}D^{\alpha}_{x}w(0,t)\\
w(x,0)=\hat{P}_{0}(x)\\
w_{x}(0,t)=w(1,t)=0\\
\end{array}\right.
\end{equation}
}
\subsection{Infinite-dimensional backstepping-like transformation and observer gains selection }\leavevmode\par\noindent
 \R{we propose to extend the backstepping approach to the space FPDE in \eqref{eq2A} for the determination of the space dependent gains $H_{1}(x)$ and $H_{2}(x)$. We start by introducing the Voltera coordinates transformation \cite{zhou2018boundary}, \cite{krstic2008boundary} and \cite{smyshlyaev2005backstepping}:
\begin{equation}\label{eq25B} 
Z(x,t)=\mathcal{V}\{\tilde{P}(.,t)\}(x):=w(x,t)-\int_{0}^{x}R(x,y)w(y,t)dy.
\end{equation}
This transformation is invertible, the formula for the inverse can be found in \cite{krstic2008boundary}. We consider the following space FPDE  target system:
\begin{equation}\label{eq26B} 
\left\{\begin{array}{llll}    
 \dfrac{\partial Z(x,t)}{\partial t}  - \frac{k(t)}{C \mu} \dfrac{\partial}{\partial x} _{}^{C}D^{\alpha}_{x}  Z(x,t)=0, \\
 Z(x,0)=\mathcal{V}\{\hat{P}_{0}(x)\}(x):=Z_{0}(x)
\\
Z(1, t) =Z_{x}(0, t)=0,   \\
\end{array}\right.
\end{equation}
and the kernel system:
\begin{equation}\label{eq27B} 
\left\{\begin{array}{llll}    
\dfrac{\partial}{\partial x} _{y}^{C}D^{\alpha}_{x} R(x,y)=_{}^{C}D^{\alpha}_{y,x} \dfrac{\partial}{\partial y} R(x,y), \\
_{y}^{C}D^{\alpha-1}_{x} R(x,y)|_{y=x}=0,\\
_{}^{C}D^{\alpha}_{y,x} R(x,y)|_{y=x}=0,\\
R(x,x)=0,   \\
\end{array}\right.
\end{equation}
where, ${}^{C}D^{\alpha}_{y,x}$ is the right sided Caputo derivative defined in \eqref{sided}. With the extra conditions that allow the choice of the observer gains:
\begin{equation}\label{eq27BB}
\left\{\begin{array}{llll}    
_{}^{C}D^{\alpha}_{y,x} R(x,y)|_{y=0}=\mathcal{V}\{H_{1}\}(x),\\
 R(x,y)|_{y=0}=-\mathcal{V}\{H_{2}\}(x),\\
\end{array}\right.
\end{equation}
\begin{theorem}
Using transformation \eqref{eq25B}, if there exist a twice continuously differentiable kernel function $R(x, y)$ satisfying \eqref{eq27B}  then, system \eqref{eq2A} is equivalent to \eqref{eq26B}.
\end{theorem}}
\begin{proof}

We start by computing the time classical derivative of the new coordinate $Z(x, t)$ using the same derivations as in \eqref{eq28}:
\begin{equation}\label{eq28B} 
\begin{split}
\frac{\partial}{\partial t} Z(x, t)=  \frac{\partial}{\partial t} w(x,t)-\int_{0}^{x}R(x,y) \frac{\partial}{\partial t} w(y,t)dy \\
=\frac {k(t)}{C \mu} \frac{\partial}{\partial x} _{}^{C}D^{\alpha}_{x}w(x,t)+w(0,t) \mathcal{V} \{H_{1}\}(x)+_{}^{C}D^{\alpha}_{x}w(0,t) \mathcal{V} \{H_{2}\}(x)  \\\
 - \frac{k(t)}{C \mu}  \int_{0}^{x}  \hspace{1mm} _{}^{C}D^{\alpha}_{y,x}  \frac{\partial}{\partial y} R(x,y) w(y,t) dy  \quad \\
+ \frac{k(t)}{C \mu}  [    _{}^{C}D^{\alpha-1}_{y,x}   \frac{\partial}{\partial y} R(x,y)  w(y,t) ]_{y=0}^{y=x}\\
- \frac{k(t)}{C \mu} [R(x,y) _{}^{C}D^{\alpha}_{y} w(y,t)]_{y=0}^{y=x}\\
 \end{split}                                    
\end{equation}
where, \eqref{eq28B} is obtained by first applying a classical integration by parts then the fractional integration by parts \eqref{eq2}. 
Using the generalized Leibnitz differentiation rule \eqref{eq3} and the same derivations as in \eqref{eq29} , we obtain the Caputo spatial fractional derivative of the new coordinate $Z(x, t)$:
\begin{equation}\label{eq29B} 
\begin{split}
_{}^{C}D^{\alpha}_{x}  Z(x, t)= _{}^{C}D^{\alpha}_{x} w(x,t) \\
-\int_{0}^{x} \hspace{1mm} _{y}^{C}D^{\alpha}_{x} R(x,y) w(y,t)dy\\
-[_{y}^{C}D^{\alpha-1}_{x} R(x,y)w(y, t)]|_{y=x}, \\
\end{split}
\end{equation}
then, by the classical Leibnitz differentiation rule, we obtain:
\vspace{-3mm}
\begin{equation}\label{eq30B} 
\begin{split}
\frac{\partial}{\partial x}   _{}^{C}D^{\alpha}_{x} Z(x, t)= \frac{\partial}{\partial x}  _{}^{C}D^{\alpha}_{x} w(x,t)  \\
-\int_{0}^{x} \frac{\partial}{\partial x}   _{y}^{C}D^{\alpha}_{x} R(x,y) w(y,t)dy\\
- [_{y}^{C}D^{\alpha}_{x} R(x,y)w(y, t)]|_{y=x}\\
-\frac{\partial}{\partial x} ( [_{y}^{C}D^{\alpha-1}_{x} R(x,y)w(y, t)]|_{y=x}).\\
\end{split}
\end{equation}
Thus, by \eqref{eq4}, \eqref{eq25B}, \eqref{eq27B}, \eqref{eq28B}, \eqref{eq29B} and \eqref{eq30B} we have:
\begin{equation}\label{eq31B} 
\frac{\partial Z(x,t)}{\partial t} \frac{k(t)}{C \mu} \frac{\partial}{\partial x} _{}^{C}D^{\alpha}_{x} Z(x,t)=0
\end{equation}
Which completes the proof.
\end{proof}

\R{
\begin{remark}
The kernel system \eqref{eq27B} admits at least one family of twice continuously differentiable solutions $K(x,y)$ in the triangle $ 0 \leq y \leq x \leq 1$ and which is given by \eqref{eq33}.
\end{remark}
From the conditions in \eqref{eq27BB} we have that, the gains $H_{1}(x)$ and $H_{2}(x)$ should be chosen as follows: 
$$
H_{1}(x)=\mathcal{V}^{-1}\{_{}^{C}D^{\alpha}_{y,x} R(x,y)|_{y=0}\}
$$ 
and 
$$
H_{2}(x)=-\mathcal{V}^{-1}\{ R(x,y)|_{y=0}\}
$$
\begin{theorem}\label{obs}
 System \eqref{eq27B} admits a unique solution given by:
\begin{equation*} 
Z(x, t)=\int _{- \infty}^{+ \infty}  \mathcal{G}_{\alpha}(\mid x-y \mid,t-\tau)  Z_{0}(y)dy 
\end{equation*}
 if $Z_{0}(x) \in L^{1}(\mathbb{R})$ (which means that $\hat{P}_{0}(x) \in L^{1}(\mathbb{R})$ ) and $k(t)$ is bounded from bellow $(k_{0} \leq k(t))$. Then, the observer target system given in \eqref{eq27B} is asymptotically stable in $\B{L_{\infty}(\mathbb{R})}$. Furthermore, $w(x,t)$ is asymptotically stable in $\B{L_{\infty}(\mathbb{R})}$as well:
\begin{equation*} 
\lim_{t \rightarrow +\infty} \norm{w(x, t)}_{\infty}=0 
\end{equation*}
\end{theorem}
\begin{proof}
Using \R{T}heorem \ref{thm1}, system \eqref{eq27B} admits a unique solution given by:
\begin{equation*} 
Z(x, t)=\int _{- \infty}^{+ \infty}  \mathcal{G}_{\alpha}(\mid x-y \mid,t-\tau)  Z_{0}(y)dy 
\end{equation*}
Using \R{T}heorem \ref{thm2}, because $Z_{0}(x) \in L^{1}(\mathbb{R})$ and $k(t)$ is bounded from bellow $(k_{0} \leq k(t))$. Then, the observer target system given in \eqref{eq27B} is asymptotically stable. Thus:
\begin{equation} \label{00}
\lim_{t \rightarrow +\infty} \norm{Z(x, t)}_{\infty}=0 
\end{equation}
Thus using the inverse of the transformation \eqref{eq25B} given in [] by:
\begin{equation}\label{inv}
w(x,t)=Z(x,t)+\int_{0}^{x}\bar{R}(x,y)Z(y,t)dy
\end{equation}
taking the sup over $[0,x]$ of the absolute value of \eqref{inv} and using Schwarz inequality, we get:
\begin{equation}\label{invv}
\norm{w(x,t)}_{\infty}\leq\norm{Z(x,t)}_{\infty}+\norm{\bar{R}(x,y)}_{2}\norm{Z(y,t)}_{2}
\end{equation}
Then, we take the limit of \eqref{invv} as $t \rightarrow +\infty$:
\begin{equation}\label{invvv}
\lim_{t \rightarrow +\infty}\norm{w(x,t)}_{\infty}\leq \lim_{t \rightarrow +\infty} \norm{Z(x,t)}_{\infty}+\norm{\bar{R}(x,y)}_{2} \lim_{t \rightarrow +\infty}\norm{Z(y,t)}_{2}
\end{equation}
using \eqref{00} and the equivalence of the $L^{2}$ and the $L^{\infty}$ norms we get that:
\begin{equation*}
\lim_{t \rightarrow +\infty}\norm{w(x,t)}_{\infty}=0
\end{equation*}
\end{proof}
The result of Theorem \ref{obs} is quite interesting but, we still need to prove that the observer error $\tilde{P}(x,t)$ is asymptotically stable as well. Actually, in view of \eqref{like}, one has to show that $\tilde{\theta}(t)$ is also exponentially vanishing and $|\Lambda(x,t)|$ is bounded. Before that, we will first investigate the selection of the estimate of the disturbances $\hat{\theta}(t)$}
\subsection{Disturbances adaptive law selection}
\R{The choice of the estimate of the disturbances is model free thanks to the choice of $\Lambda(x,t)$ which allowed the rejection of the effect of the disturbances from the state estimation process. That is why, we propose the parameter adaptive law in \cite{ahmed2016adaptive} which is enough to guarantee the exponential stability of the disturbances estimation error $\tilde{\theta}(t)$ independently of the state estimation problem. 
\begin{theorem}\label{0000}
Using the following parameter adaptive law parameter adaptive law in \cite{ahmed2016adaptive}:
\begin{equation*}
\begin{split}
\dot{\hat{\theta}}(t)=\frac{R(t)\Lambda^{T}(0,t)}{1+\Lambda(0,t)^{T}\Lambda(0,t)}\tilde{P}(0,t)\\
\dot{R}(t)=R(t)-\frac{R(t)\Lambda^{T}(0,t)\Lambda(0,t)R(t)}{1+\Lambda(0,t)^{T}\Lambda(0,t)}
\end{split}
\end{equation*}
where, $R(t) \in \mathbb{R}^{2\times2}, \quad  \hat{\theta}(0)$ \text{and}  $R(0)=R_{0}$ are arbitrarily chosen with $R_{0}=R_{0}^{T}>0$. 
the disturbance estimation error $\tilde{{\theta}}(t)$ is exponentially stable. 
\end{theorem}
\begin{proof}
consider the following parameter adaptive law parameter adaptive law in \cite{ahmed2016adaptive}:
\begin{equation}\label{law}
\begin{split}
\dot{\hat{\theta}}(t)=\frac{R(t)\Lambda^{T}(0,t)}{1+\Lambda(0,t)^{T}\Lambda(0,t)}\tilde{P}(0,t)\\
\dot{R}(t)=R(t)-\frac{R(t)\Lambda^{T}(0,t)\Lambda(0,t)R(t)}{1+\Lambda(0,t)^{T}\Lambda(0,t)}
\end{split}
\end{equation}
The choice of using only the measurements $z_{m}(t)$ for the disturbances estimation will prove its efficiency later. This parameter adaptive law is a variant of the least squares estimator, commonly referred to forgetting factor least squares \cite{ioannou2012robust}. Thus \eqref{law} is equivalent to:
\begin{equation}\label{laww}
\frac{d R^{-1}}{dt}=-R^{-1}(t)+\frac{\Lambda^{T}(0,t)\Lambda(0,t)}{1+\Lambda(0,t)^{T}\Lambda(0,t)}
\end{equation}
 using transformation \eqref{like}, \eqref{law} equivalent to:
 \begin{equation}\label{lawww}
\dot{\tilde{\theta}}(t)=-\frac{R(t)\Lambda^{T}(0,t)\Lambda(0,t)}{1+\Lambda(0,t)^{T}\Lambda(0,t)}\tilde{\theta}(t)+\frac{R(t)\Lambda^{T}(0,t)}{1+\Lambda(0,t)^{T}\Lambda(0,t)}w(0,t)
\end{equation}
Using the following Lyapunov function:
 \begin{equation}\label{lya}
V_{1}(t)=\tilde{\theta}^{T}(t)R^{-1}(t) \tilde{\theta}(t)
\end{equation}
using \eqref{law},\eqref{laww} and \eqref{lawww}, it has been proved in \cite{ahmed2016adaptive} that:
 \begin{equation}\label{lyapunov}
\dot{V}_{1}(t)=\tilde{\theta}^{T}(t)\dot{R}^{-1}(t) \tilde{\theta}(t)+2\tilde{\theta}^{T}(t)R^{-1}(t) \dot{\tilde{\theta}}(t)\leq - V_{1}(t)+w^{2}(0,t)
\end{equation}
using Theorem \ref{obs}, $w^{2}(0,t)$ is asymptotically vanishing. Thus, from \eqref{lyapunov} and by the comparison lemma \cite{khalil2002nonlinear}, $V_{1}(t)$ is exponentially vanishing. In view of \eqref{lya} so is $\tilde{\theta}(t)$.
\end{proof}
\begin{theorem}
the observer error $\tilde{P}(x,t)$ in \eqref{eq2A} is asymptotically stable in $\B{L_{\infty}(\mathbb{R})}$. That is,
\begin{equation*} 
\lim_{t \rightarrow +\infty} \norm{\tilde{P}(x, t)}_{\infty}=0 
\end{equation*}
\end{theorem}
\begin{proof}
We start by the transformation in \eqref{like}
\begin{equation}\label{like1} 
\tilde{P}(x,t)=w(x,t)+\Lambda(x,t)\tilde{\theta}(t),
\end{equation}
taking the sup over $[0,x]$ of the absolute value of \eqref{like1} and using Cauchy-schwarz inequality, we get:
\begin{equation}\label{like2} 
\norm{\tilde{P}(x,t)}_{\infty}\leq\norm{w(x,t)}_{\infty}+|\tilde{\theta}(t)|\norm{\Lambda(x,t)}_{\infty},
\end{equation}
by taking the limit of \eqref{like2} as $t \rightarrow +\infty$ and using the fact that $\Lambda(x,t)$ is stable (Theorem \ref{00000}) and that $\tilde{\theta}(t)$ vanishes exponentially (Theorem \ref{0000})we get the result 
\end{proof}}

\section{Observer Based Output Regulation}
Since we obtain an approximated state $\hat{P}(x, t)$ from the output by observer \eqref{eq1A}, it follows from the output controllers \eqref{eq40} and \eqref{eq40N} that an observer-based output regulation controller should be designed as:
\begin{equation}\label{co1}
u(t)=P(1, t)=\int_{0}^{1}k(1,y)\hat{P}(y,t)dy+[M^{T}(x)]|_{x=1}V(t).
\end{equation}
and 
\begin{equation}\label{co2}
u(t)=P(1, t)=\int_{0}^{t}\hat{P}(1,\tau)l(\tau,t)d\tau+[M^{T}(x)]|_{x=1}V(t).
\end{equation}
respectively. 
\begin{theorem}
Under feedback controllers \eqref{co1} (respectively \eqref{co2}), \R{if $\hat{P}_{0}(x) \in L^{1}(\mathbb{R})$ and $k(t)$ is bounded from bellow $(k_{0} \leq k(t))$}, we have that the closed loop system: 
\R{
 \begin{equation}\label{eq1C}
\left\{\begin{array}{llll} 
\frac{\partial P(x,t)}{\partial t}  - \frac{k(t)}{C \mu} \frac{\partial}{\partial x} _{}^{C}D^{\alpha}_{x}  P(x,t)=\frac{1}{C \rho}d_{1}(t)f(x),\\
 P(x, 0)=0,
\\
P(1, t) = \int_{0}^{1}K(1,y)\hat{P}(y,t)dy+[M^{T}(x)]|_{x=1}V(t),  \\  P_{x}(0, t) = d_{2}(t),\\ P(0,t)=z_{m}(t), _{}^{C}D^{\alpha}_{x} P(0,t)=y_{m}(t), \\
 \frac{\partial \hat{P}(x,t)}{\partial t}  - \frac{k(t)}{C \mu} \frac{\partial}{\partial x} _{}^{C}D^{\alpha}_{x}  \hat{P}(x,t)=\frac{1}{C \rho}\hat{d}_{1}(t)f(x)+u(x,t)\\
 \qquad \qquad -\frac{k(t)}{C \mu}H_{1}(x)(\hat{P}(0,t)-z_{m}(t))-\frac{k(t)}{C \mu}H_{2}(x)(_{}^{C}D^{\alpha}_{x} \hat{P}(0,t)-y_{m}(t)),\\
 \hat{P}(x, 0)=\hat{P}_{0}(x),\\
 \hat{P}(1, t) = \int_{0}^{1}K(1,y)\hat{P}(y,t)dy+[M^{T}(x)]|_{x=1}V(t),  \\
  \hat{P}_{x}(0, t) = \hat{d}_{2}(t),   \\
\end{array}\right.
\end{equation}
is asymptotically stable. \R{Where, $\hat{d}_{1}(t)$ and $\hat{d}_{2}(t)$ have to satisfy the adaptive law in \eqref{law}}}
\end{theorem}
\begin{proof}
\R{Using the observer error variable  $\tilde{P}=\hat{P}(x,t)-P(x,t)$, \eqref{eq1C} becomes:
 \begin{equation}\label{eq1C*}
\left\{\begin{array}{llll} 
\frac{\partial P(x,t)}{\partial t}  - \frac{k(t)}{C \mu} \frac{\partial}{\partial x} _{}^{C}D^{\alpha}_{x}  P(x,t)=\frac{1}{C \rho}d_{1}(t)f(x),\\
 P(x, 0)=0,
\\
P(1, t) = \int_{0}^{1}K(1,y)\hat{P}(y,t)dy+[M^{T}(x)]|_{x=1}V(t),  \\  P_{x}(0, t) = d_{2}(t),\\ P(0,t)=z_{m}(t), _{}^{C}D^{\alpha}_{x} P(0,t)=y_{m}(t), \\
 \frac{\partial \tilde{P}(x,t)}{\partial t}   - \frac{k(t)}{C \mu} \frac{\partial}{\partial x} _{}^{C}D^{\alpha}_{x}  \tilde{P}(x,t)=\frac{1}{C \rho}\tilde{d}_{1}(t)f(x)+u(x,t)\\
 \qquad \qquad -\frac{k(t)}{C \mu}H_{1}(x)\tilde{P}(0,t)-\frac{k(t)}{C \mu}H_{2}(x)_{}^{C}D^{\alpha}_{x}\tilde{P}(0,t) ,\\
 \tilde{P}(x, 0)=\hat{P}_{0}(x), \tilde{P}(1, t) =0,  \\
  \tilde{P}_{x}(0, t) =\tilde{d}_{2}(t),   \\
\end{array}\right.
\end{equation}
Using the transformations in  \eqref{eq25} and \eqref{eq36},  system \eqref{eq1C*} is equivalent to:
 \begin{equation}\label{eq2Cx}
\left\{\begin{array}{llll} 
\frac{\partial e(x,t)}{\partial t}  = \frac{k(t)}{C \mu} \frac{\partial}{\partial x} _{}^{C}D^{\alpha}_{x} e(x,t) \\
 e(x, 0) =e_{0}(x), \qquad e_{x}(0, t) =0, \\
 e(1, t) =\int^{1}_{0}K(1,y)[\mathcal{V}^{-1}\{Z(y,t)\}]dy,   \\
  \frac{\partial \tilde{P}(x,t)}{\partial t}   - \frac{k(t)}{C \mu} \frac{\partial}{\partial x} _{}^{C}D^{\alpha}_{x}  \tilde{P}(x,t)=\frac{1}{C \rho}\tilde{d}_{1}(t)f(x)+u(x,t)\\
 \qquad \qquad -\frac{k(t)}{C \mu}H_{1}(x)\tilde{P}(0,t)-\frac{k(t)}{C \mu}H_{2}(x)_{}^{C}D^{\alpha}_{x}\tilde{P}(0,t) ,\\
 \tilde{P}(x, 0)=\hat{P}_{0}(x), \tilde{P}(1, t) =0,  \\
  \tilde{P}_{x}(0, t) =\tilde{d}_{2}(t),   \\
\end{array}\right.
\end{equation}
where, $[\mathcal{V}^{-1}\{Z(y,t)\}] =\tilde{P}(x,t).$. Using the transformations in \eqref{like} and \eqref{eq25B}, system \eqref{eq2Cx} is equivalent to:
 \begin{equation}\label{eq2C}
\left\{\begin{array}{llll} 
\frac{\partial e(x,t)}{\partial t}  = \frac{k(t)}{C \mu} \frac{\partial}{\partial x} _{}^{C}D^{\alpha}_{x} e(x,t) \\
 e(x, 0) =e_{0}(x), \qquad e_{x}(0, t) =0, \\
 e(1, t) =\int^{1}_{0}K(1,y)[\mathcal{V}^{-1}\{Z(y,t)\}]dy,   \\
  \frac{\partial Z(x,t)}{\partial t}  - \frac{k(t)}{C \mu} \frac{\partial}{\partial x} _{}^{C}D^{\alpha}_{x}  Z(x,t)=0, \\
 Z(x,0)=Z_{0}(x)
\\
 Z(1, t) =Z_{x}(0, t) =0 \\
\end{array}\right.
\end{equation}}
 \R{It is clear from \eqref{eq2C} that the observer error $Z(x,t)$ is asymptotically stable using Theorem \ref{obs} (because,  $\hat{P}_{0}(x) \in L^{1}(\mathbb{R})$ and $k(t)$ is bounded from bellow $(k_{0} \leq k(t))$). Thus, the output tracking error $e(x,t)$ is asymptotically stable as well using Theorem \ref{error} (no conditions required on $e_{0}(x)$)}.
\end{proof}
\section{\R{Numerical results}}
\begin{figure}[!ht]
	\centering
	\begin{subfigure}[b]{0.49\textwidth}
         \centering
         \includegraphics[width=\textwidth]{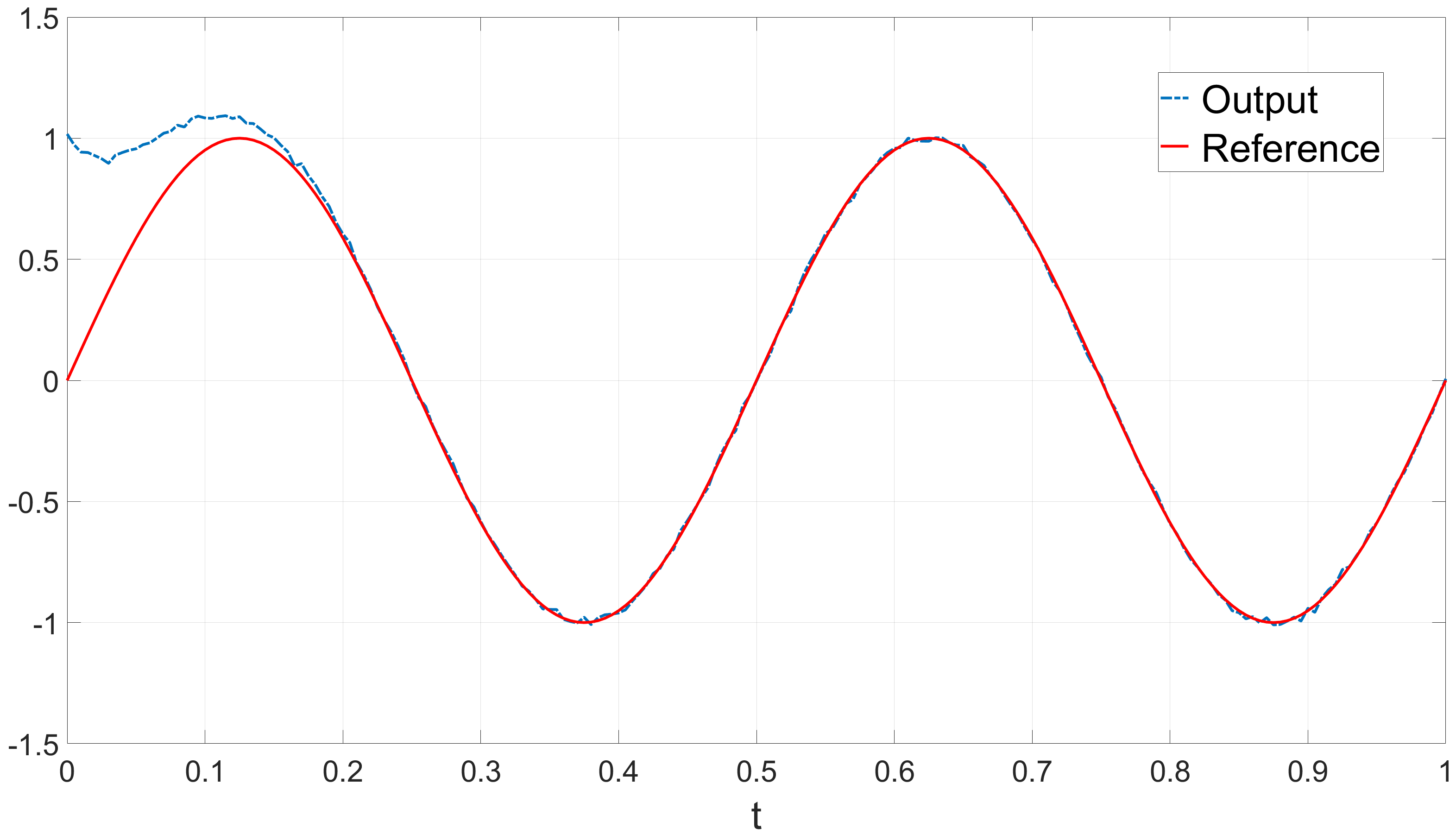}
     \end{subfigure}
     \hfill
     \begin{subfigure}[b]{0.49\textwidth}
         \centering
         \includegraphics[width=\textwidth]{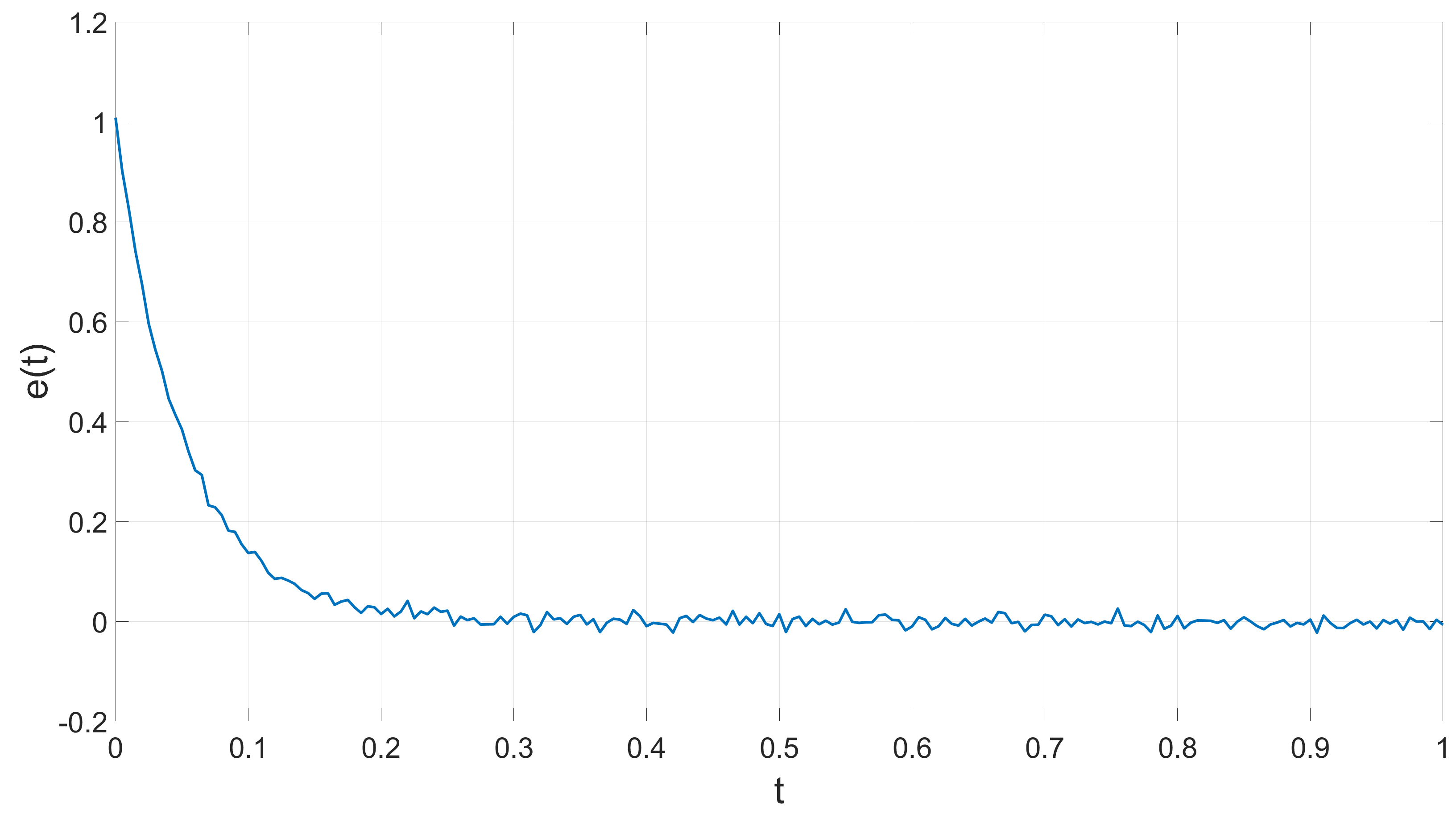}
     \end{subfigure}
     \caption{\R{Reference tracking process (left) and reference tracking error (right) for $\sigma=0.01$.} }\label{fig12}  
\end{figure}

\begin{figure}[!ht]
	\centering
	\begin{subfigure}[b]{0.49\textwidth}
         \centering
         \includegraphics[width=\textwidth]{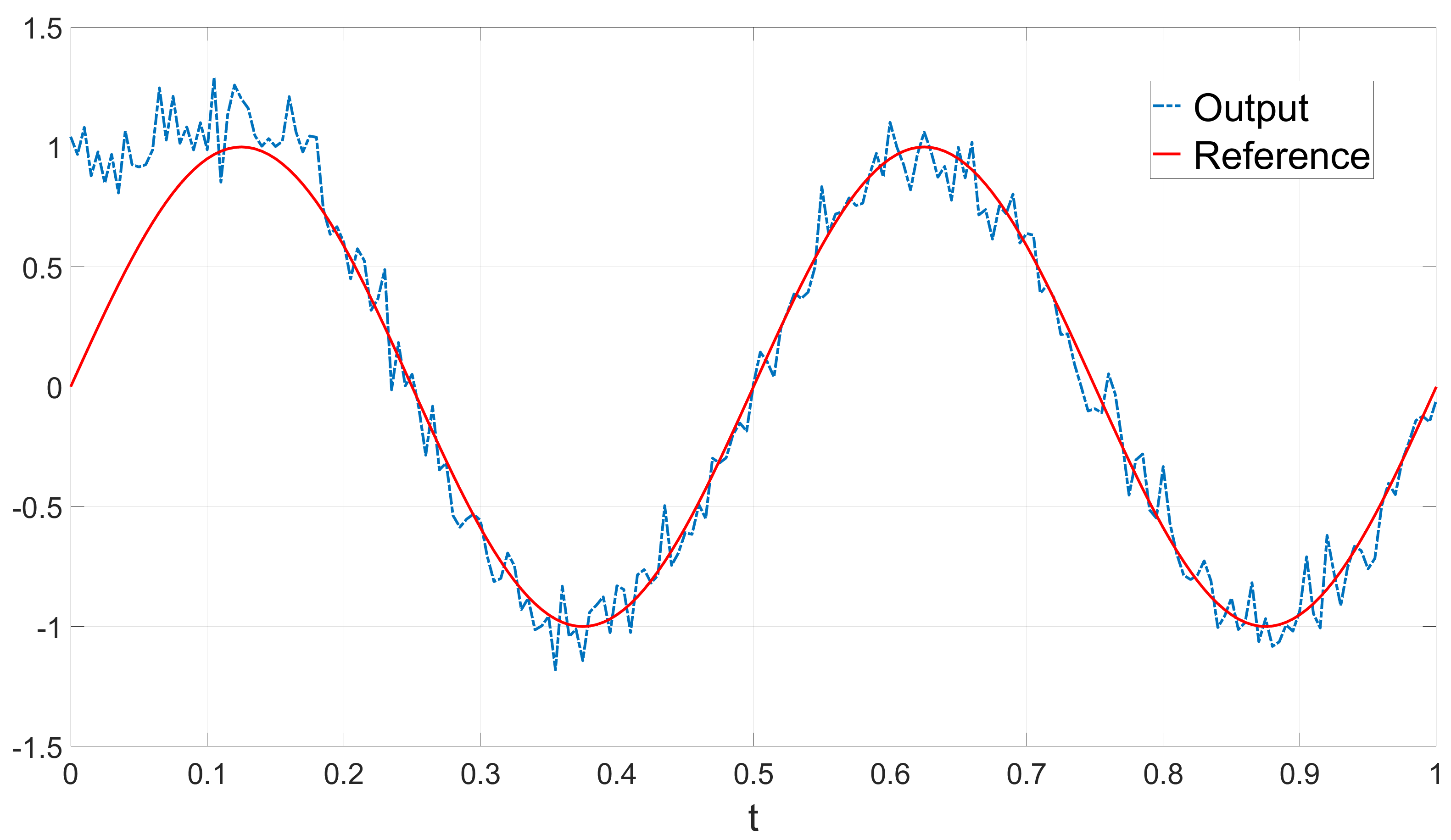}
            \end{subfigure}
     \hfill
     \begin{subfigure}[b]{0.49\textwidth}
         \centering
         \includegraphics[width=\textwidth]{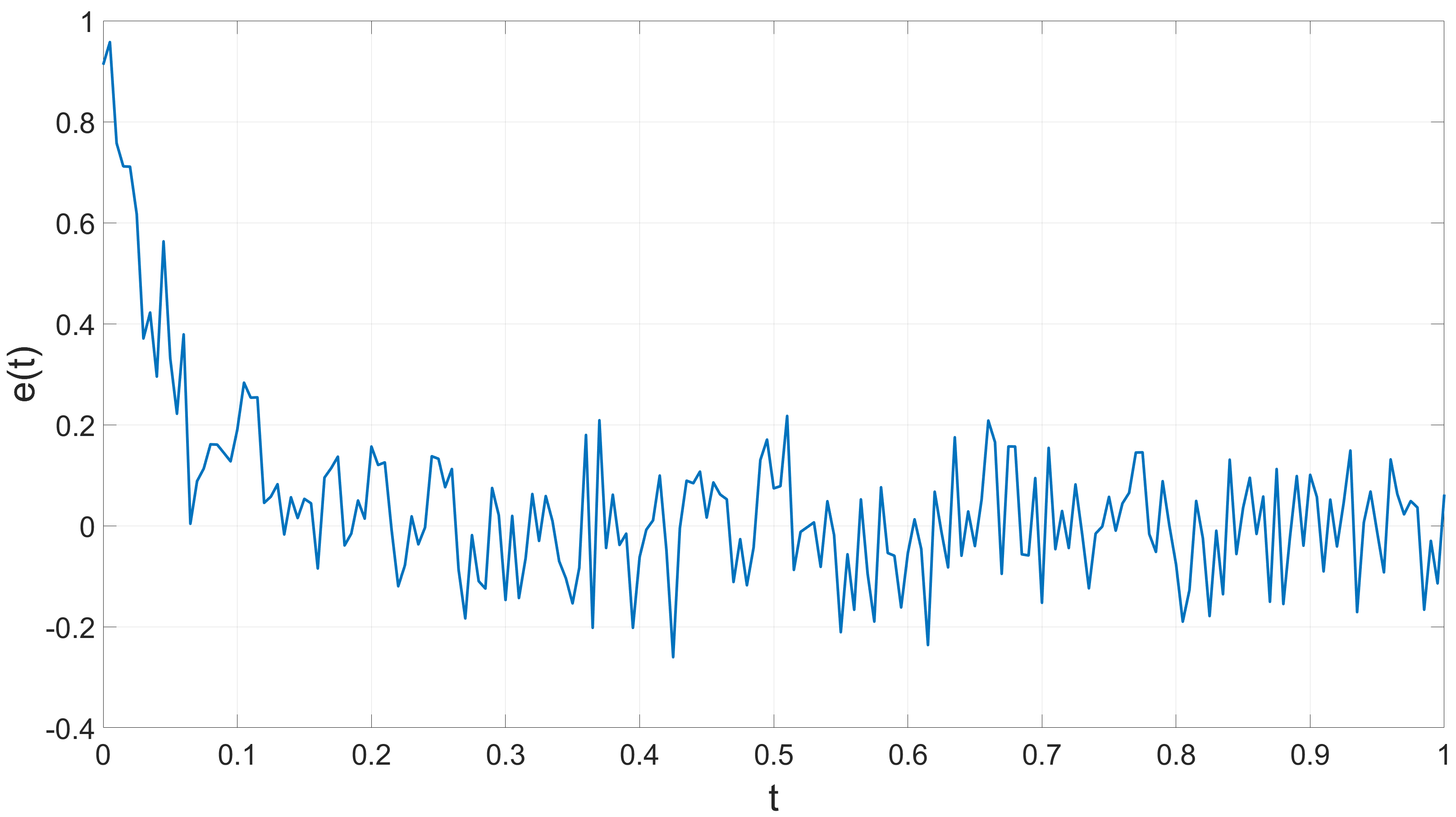}
     \end{subfigure}
     \caption{\R{Reference tracking process (left) and reference tracking error (right) for $\sigma=0.1$.} }\label{fig34}  
\end{figure}

\R{In this section, we present some numerical results to show the efficiency of the presented method to solve the reference tracking problem and for the adapted observer design.\\}
\R{we consider the following state  $P(x,t)=4e^{-t}sin(2 \pi x)$ and the following system parameters: $k=5, \mu=C=L=1, \rho=1.0726.$ Signal $V(t)$ in system \eqref{eq23} is chosen such that $\mathcal{S}=-25, V_{0}=1$. The reference is chosen as follows: $y_{d}(t)=sin(2 \pi t)$. The choice of the kernel functions $K(x,y)$ and $R(x,y)$ given by \eqref{eq33} is important since it affects the efficiency
of the algorithm. In this regard, we propose using polynomial kernel functions that satisfy systems \eqref{eq27} and \eqref{eq27B} (for $m=1$), and for which the fractional derivatives are easy to calculate, of the following form:}\vspace{-3mm}
\begin{equation*}
  \R{K(x,y)=R(x,y)=(x-y)^{3},} \vspace{-3mm}
\end{equation*}
\R{whose the fractional derivative is known analytically and given by:} \vspace{-3mm}
\begin{equation*}
  \R{_{y}^{C}D^{\alpha}_{x} K(x,y)=_{y}^{C}D^{\alpha}_{x} R(x,y)=\frac{\Gamma (4)}{\Gamma (4-\alpha)}(x-y)^{3-\alpha},}\vspace{-3mm}
\end{equation*}
\R{using Theorem \ref{thm20}. Figures \ref{fig12} and \ref{fig34} show the resulting tracking behaviour when using a high order compensator and the corresponding tracking error behavior after adding Gaussian noise to the measurements with mean equal to zero and standard deviation $\sigma$. Clearly, the estimates get very close to their true variables after a transient period. The above observations confirm the theoretical asymptotic performance described in Theorem \ref{error}.}

%  \begin{figure}[ht]
% 	\centering
% 	\includegraphics[scale=0.09]{ref0.01.png}
% 	\caption{Reference tracking process for \sigma=0.01 }\label{fig1}  
% \end{figure}
%  \begin{figure}[ht]
% 	\centering
% 	\begin{overpic}[scale=0.29]{error0.01.jpg}
% 	\end{overpic}
% 	\caption{Reference tracking error for \sigma=0.01 }\label{fig2}  
% \end{figure}

%  \begin{figure}[ht]
% 	\centering
% 	\begin{overpic}[scale=0.29]{ref0.1.jpg}
% 	\end{overpic}
% 	\caption{Reference tracking process for \sigma=0.1 }\label{fig3}  
% \end{figure}

\begin{figure}[!ht]
	\centering
	\begin{subfigure}[b]{0.49\textwidth}
         \centering
         \includegraphics[width=\textwidth]{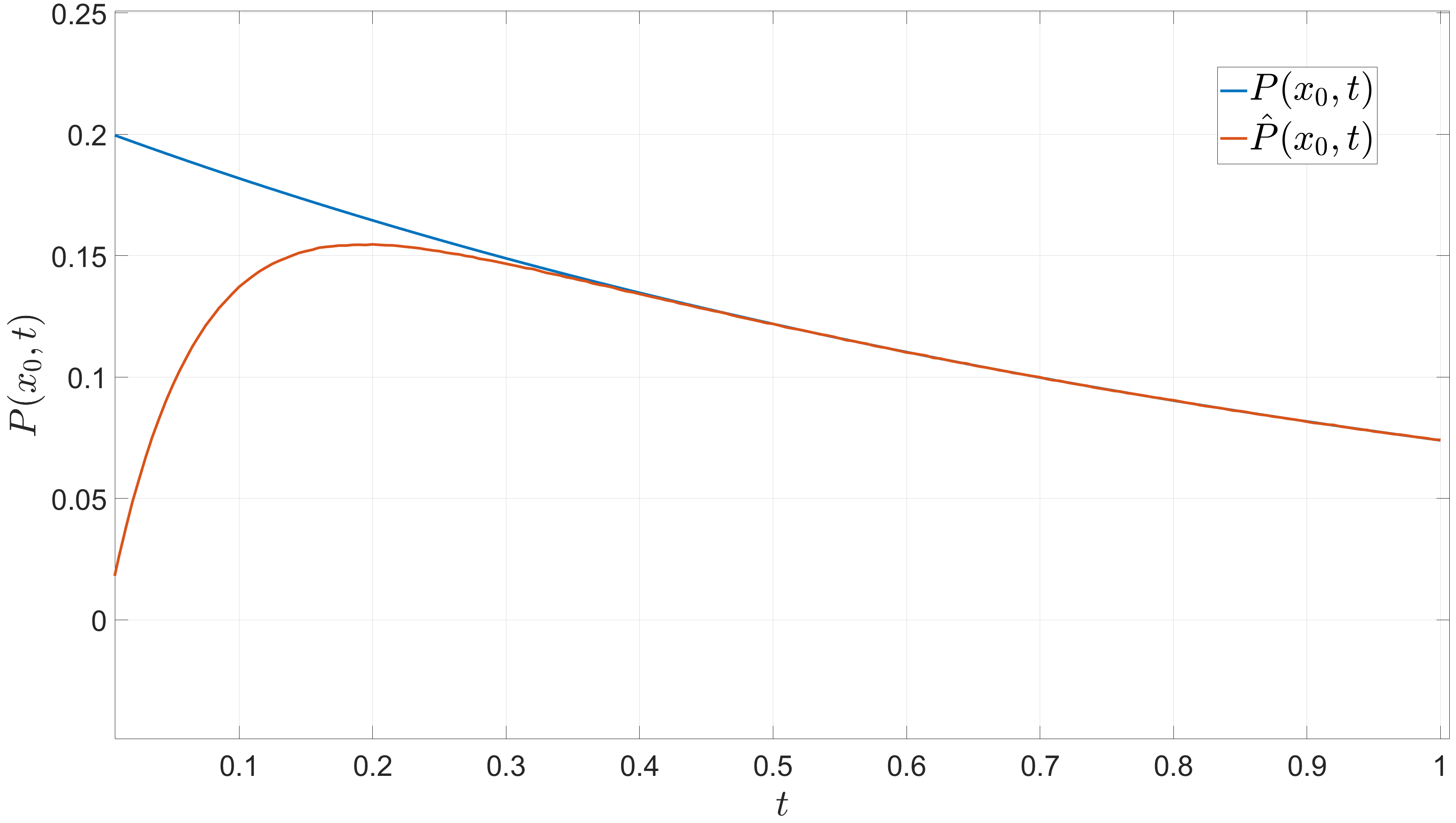}
        %  \caption{Reference tracking process} \label{fig1}
     \end{subfigure}
     \hfill
     \begin{subfigure}[b]{0.49\textwidth}
         \centering
         \includegraphics[width=\textwidth]{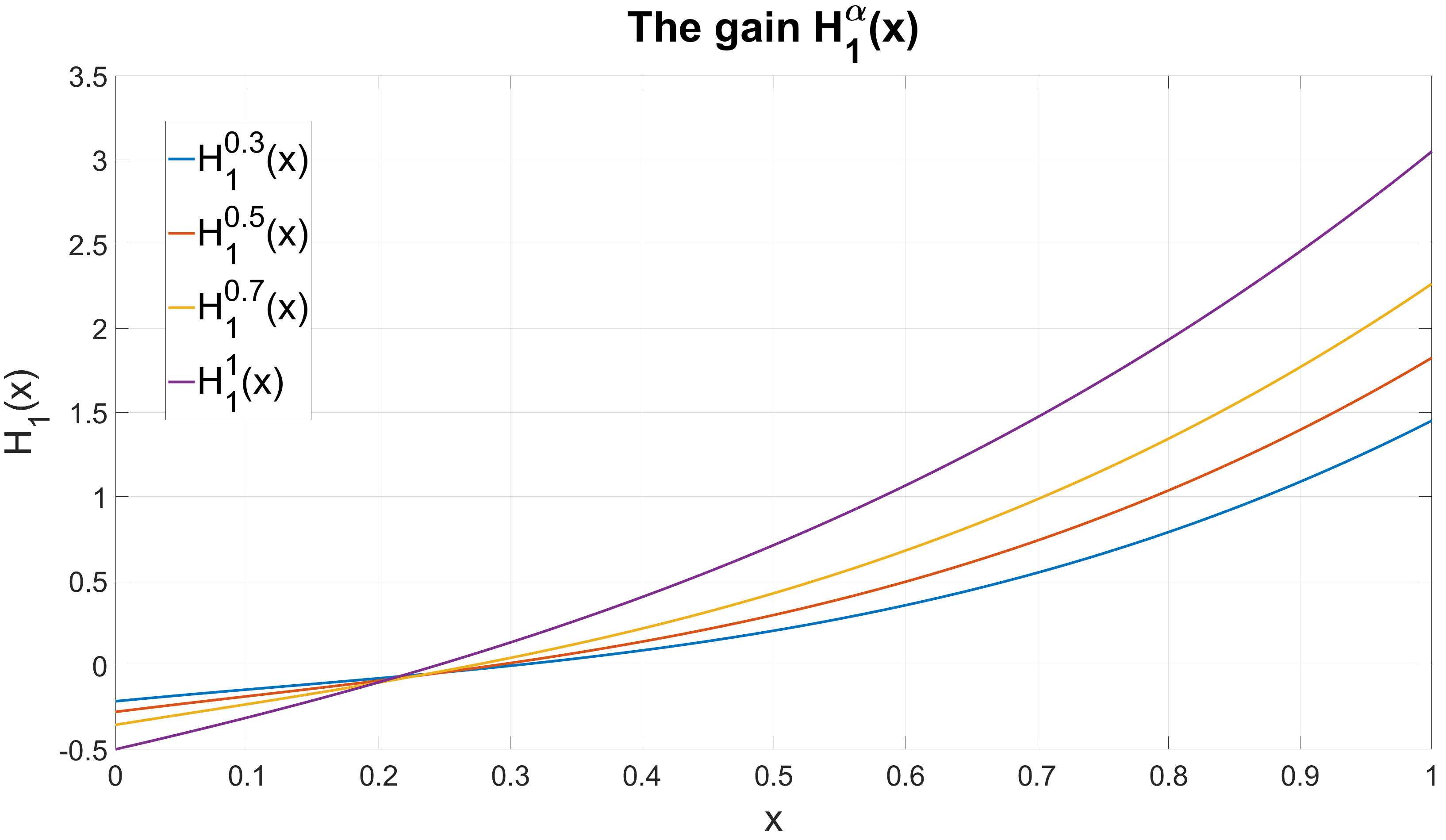}
        %  \caption{Reference tracking error} \label{fig2}
     \end{subfigure}
     \caption{\R{Time evolution of the state estimate at a particular position (left), and the gain $H_{1}(x)$ for different values of $alpha$ (right).}}\label{fig56}  
\end{figure}

\begin{figure}[!ht]
	\centering
	\begin{subfigure}[b]{0.49\textwidth}
         \centering
         \includegraphics[width=\textwidth]{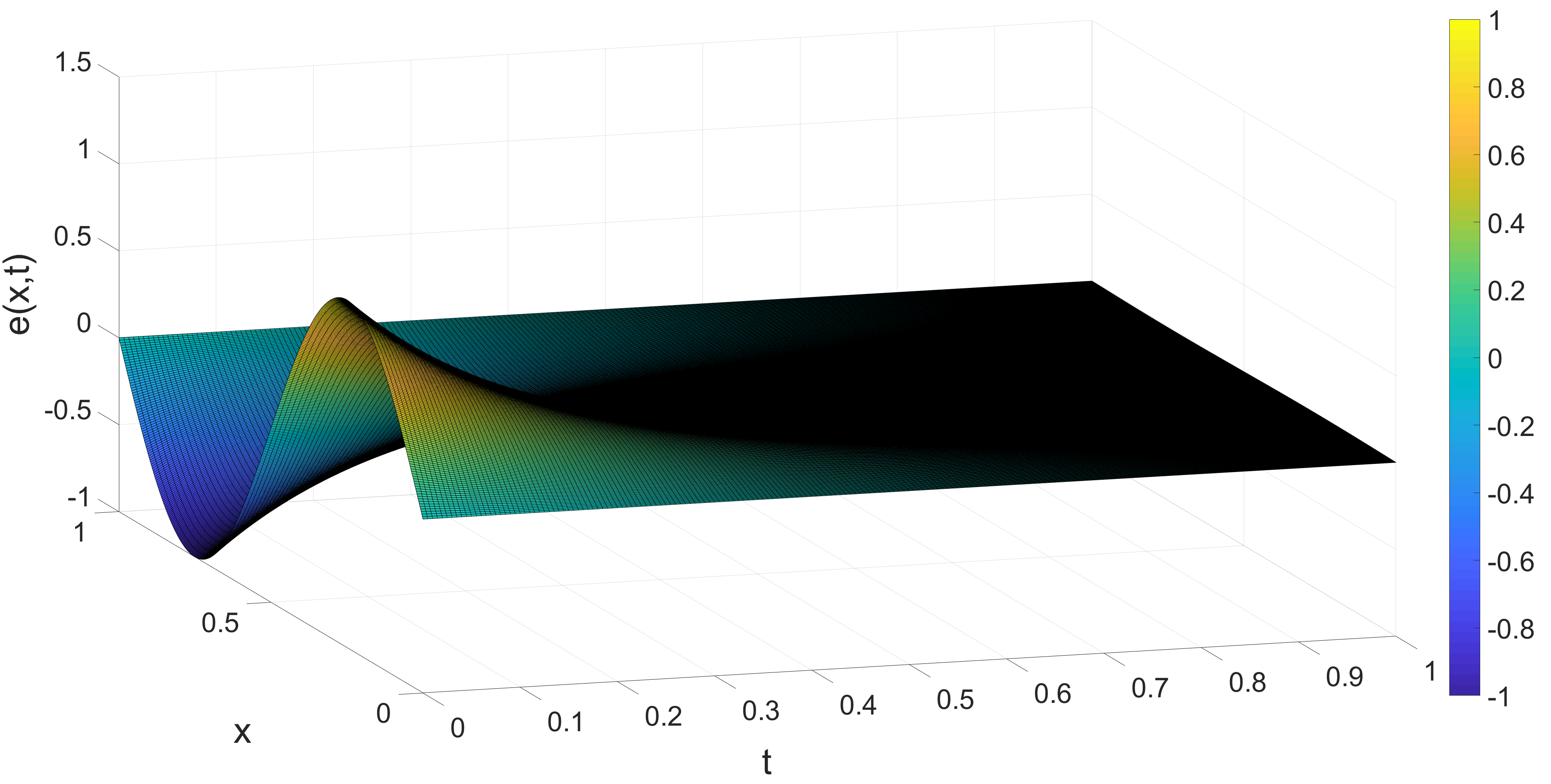}
        %  \caption{Reference tracking process} \label{fig1}
     \end{subfigure}
     \hfill
     \begin{subfigure}[b]{0.49\textwidth}
         \centering
         \includegraphics[width=\textwidth]{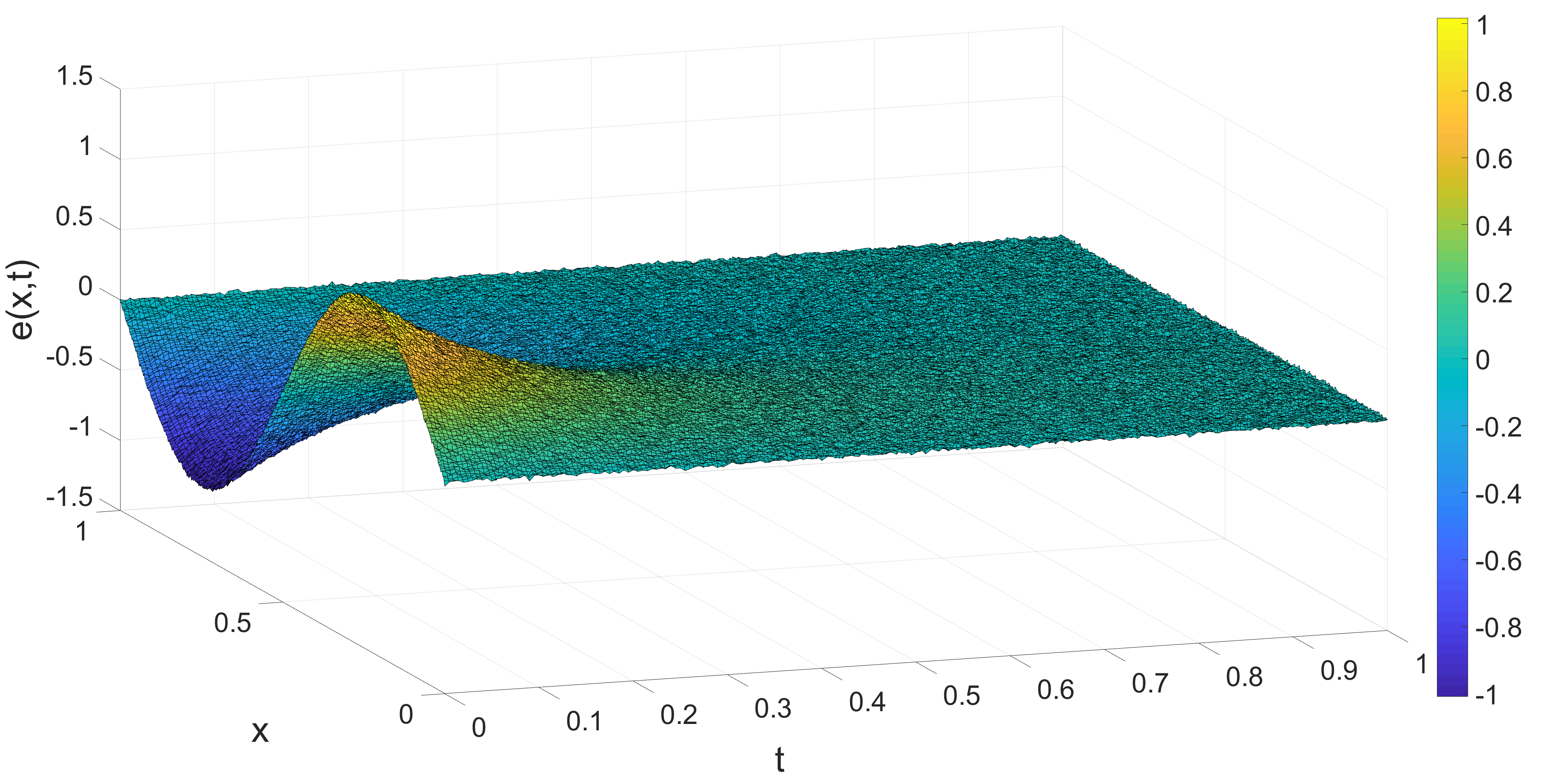}
        %  \caption{Reference tracking error} \label{fig2}
     \end{subfigure}
     \caption{\R{State estimate error and state estimate error with Gaussian noise ($\sigma=0.01).$}}\label{fig78}  
\end{figure}

\begin{figure}[!ht]
	\centering
	\begin{subfigure}[b]{0.49\textwidth}
         \centering
         \includegraphics[width=\textwidth]{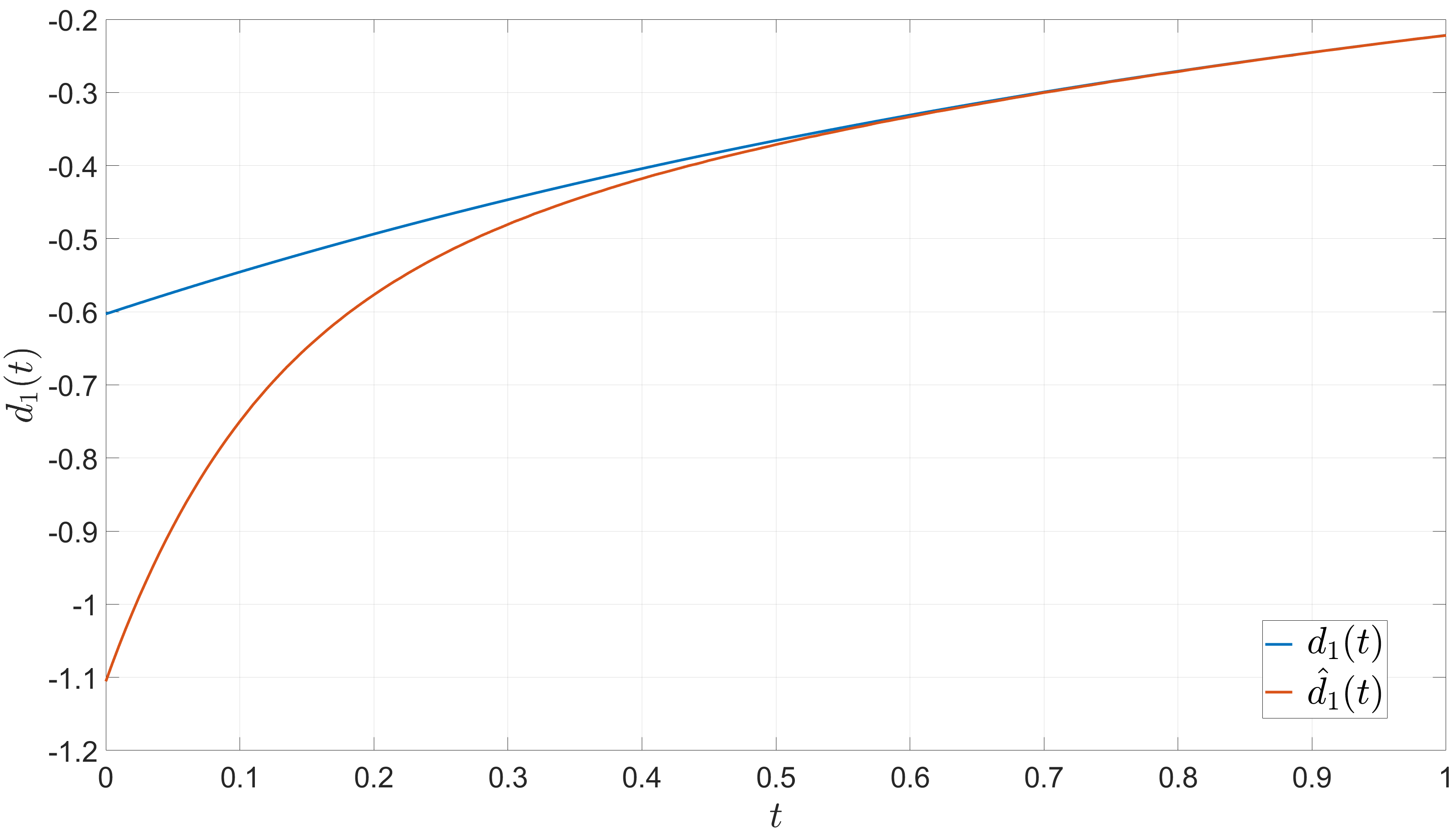}
        %  \caption{Reference tracking process} \label{fig1}
     \end{subfigure}
     \hfill
     \begin{subfigure}[b]{0.49\textwidth}
         \centering
         \includegraphics[width=\textwidth]{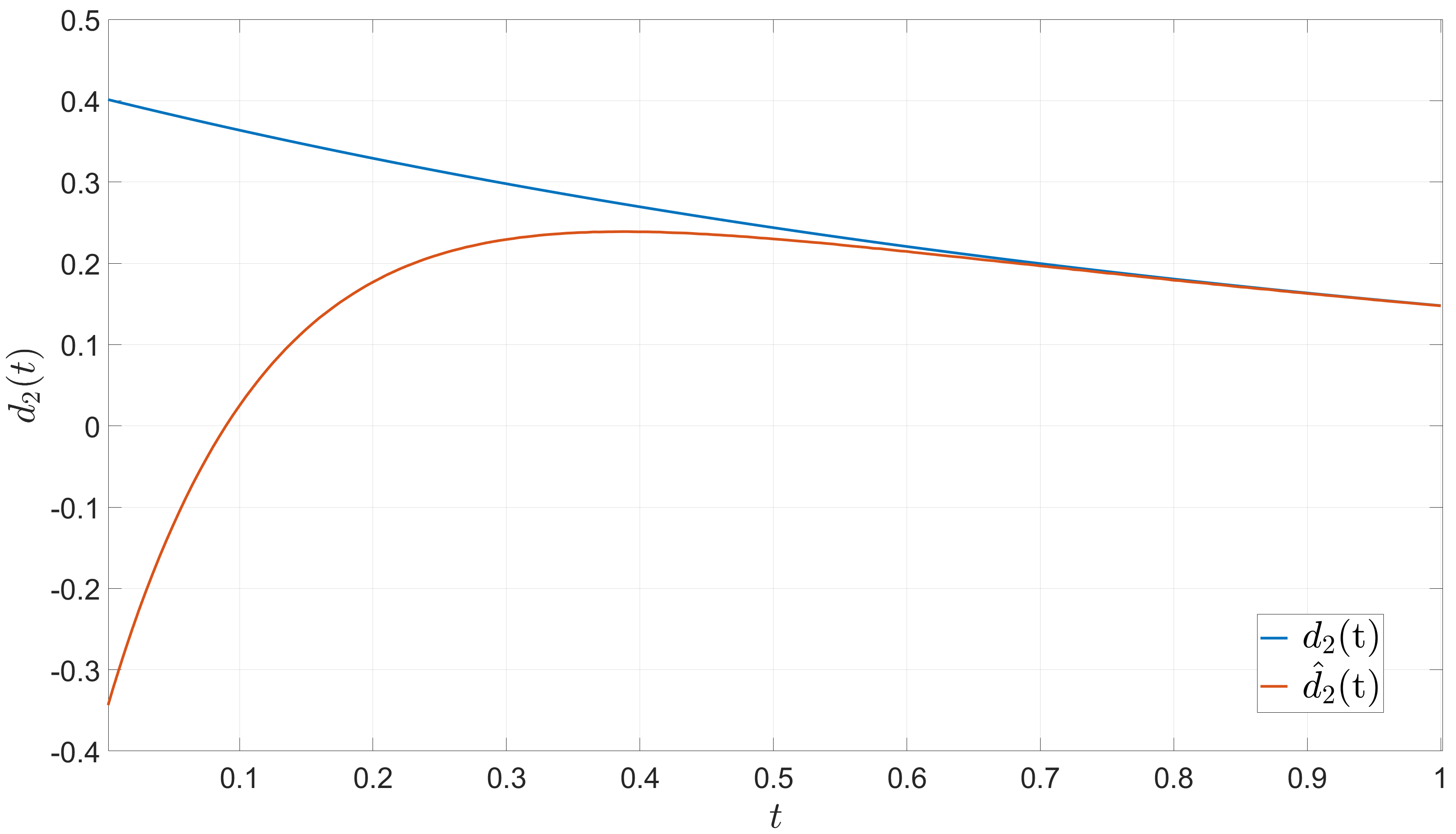}
        %  \caption{Reference tracking error} \label{fig2}
     \end{subfigure}
     \caption{\R{Estimation of the disturbances $d_{1}(t)$ and $d_{2}(t)$.}}\label{fig910}  
\end{figure}

%  \begin{figure}[ht]
% 	\centering
% 	\begin{overpic}[scale=0.29]{error0.1.jpg}
% 	\end{overpic}
% 	\caption{Reference tracking error for \sigma=0.1 }\label{fig4}  
% \end{figure}
\R{Figure \ref{fig56} (left) shows the time evolution of the state estimate at a particular position in the spatial domain, Figure \ref{fig56} (right) represents the gain $H_{1}(x)$ behavior for different values of $\alpha$. Figure \ref{fig78} shows the state estimation error for the whole space and time domains without noise (left) and with added noise (right). Figure \ref{fig910} shows the accuracy of the disturbances estimation after a transient period. The above observations confirm the theoretical asymptotic performance described in Theorems \ref{obs} and \ref{0000}. }

 %%%%%%%%%%%%%%%%%%%%%%%%%%%%%%%%%%%%%%%%% %%%%%%%%%%%%%%%%%%%%%%%%%%%%%%%%%%%%%%%%%
\section*{Conclusion}
This paper dealt with a class of boundary controlled space FPDE in the presence of disturbances describing the gas pressure in fractured media. The main contributions of this paper were first to study the stability of the considered problem. Then, to track the pressure gradient at final position in the presence of disturbances using backstepping approach. Moreover, an adaptive observer has been designed to  estimate the system's state. A fundamental solution for the considered model was given in this paper its  uniqueness has been also studied. A future direction of this work would be to   validate the obtained results by numerical simulations, and to generalize these results for the space and time FPDE. \R{A further direction could be the adaptation of the adaptive observer for a time fractional PDE and to compere to the results in \cite{ghaffour2020non} (which is a asymptotic and robust method for the state estimation for PDEs which can be applied for time fraction PDEs as well).}

\section{Appendix \R{A}} 
We recall the conventional PDE model's  derivation for the gas pressure \cite{amir2018physics}:
\begin{enumerate}
\item Mass conservation law :  
           \begin{equation*}
         \frac{\partial (\rho \varphi)(P(x,t))}{\partial t} +\nabla. (\rho u)(P(x,t))=Q(x,t),
           \end{equation*}

\item Assume the non-compressible fluid, which means that the density does not change with respect to the pressure:
           \begin{equation*}
         \rho \frac{\partial \varphi(P(x,t))}{\partial P(x,t)}  \frac{\partial P(x,t)}{\partial t}  +\rho \nabla. u(P(x,t))=Q(x,t),
           \end{equation*}
\item Assume a constant-compressible rock (a particular case of the linear-compressibility), which means that the porosity is constant with respect to the pressure:
           \begin{equation*}
         \frac{\partial P(x,t)}{\partial t}  +\frac{1}{C} \nabla. u(P(x,t))=\frac{1}{C \rho}Q(x,t), 
           \end{equation*}
\item Using the Darcy flow equation:           
          \begin{equation*}
         \frac{\partial P(x,t)}{\partial t}  +\frac{1}{C} \nabla. (-\frac{k(t)}{\mu} \frac{\partial P(x,t)}{\partial x})=\frac{1}{C \rho}Q(x,t).
           \end{equation*}
\end{enumerate}

\section{\R{Appendix B}}{\R{Proof of Theorem \eqref{00000}:}\\}
\R{Consider the  auxiliary states in \eqref{FF} with the initial and boundary conditions in \eqref{F1}:
 \begin{equation}\label{10.1}
\left\{\begin{array}{llll}   
\frac{\partial }{\partial t}\lambda_{1}(x,t)=\frac{k(t)}{C \mu} \frac{\partial}{\partial x} _{}^{C}D^{\alpha}_{x}\lambda_{1}(x,t)-H_{1}(x)\frac{\partial }{\partial x}\lambda_{1}(x,t)+\frac{1}{C \rho}f(x)\\
\lambda_{1}(x,0)=0\\
\frac{\partial }{\partial x}\lambda_{1}(x,t)|_{x=0}=\lambda_{1}(1,t)=0\\
\end{array}\right.
\end{equation}
and
\begin{equation}\label{10.2}
\left\{\begin{array}{llll}   
\frac{\partial }{\partial t}\lambda_{2}(x,t)=\frac{k(t)}{C \mu} \frac{\partial}{\partial x} _{}^{C}D^{\alpha}_{x}\lambda_{2}(x,t)-H_{1}(x)\frac{\partial }{\partial x}\lambda_{2}(x,t)\\
\lambda_{2}(x,0)=0\\
\lambda_{2}(1,t)=0\\
\frac{\partial }{\partial x}\lambda_{2}(x,t)|_{x=0}=1.\\
\end{array}\right.
\end{equation}
% We propose to solve analytically the more general system: 
% \begin{equation}\label{10.3}
% \left\{\begin{array}{llll}   
% \frac{\partial }{\partial t}\lambda(x,t)=\frac{k(t)}{C \mu} \frac{\partial}{\partial x} _{}^{C}D^{\alpha}_{x}\lambda(x,t)-H_{1}(x)\frac{\partial }{\partial x}\lambda(x,t)+F(x)\\
% \lambda(x,0)=\lambda_{0}(x)\\
% \lambda(1,t)=0\\
% \frac{\partial }{\partial x}\lambda(x,t)|_{x=0}=c.\\
% \end{array}\right.
% \end{equation}
% where, $c$ is an arbitrary constant, $F(x)$ and $\lambda_{0}(x)$ are arbitrary functions. System \eqref{10.3} is more general and cover both systems \eqref{10.1} and \eqref{10.2}. Indeed, if we take $F(x)=\frac{1}{C \rho}f(x), \lambda_{0}(x)=0$ and $c=0$, system \eqref{10.3} becomes identical to system \eqref{10.1} and if we take $F(x)=0, \lambda_{0}(x)=\lambda_{2,0}(x)$ and $c=1$, system \eqref{10.3} becomes identical to system \eqref{10.2}. So to find the analytical solutions of systems \eqref{10.1} and \eqref{10.2} it is enough to solve system \eqref{10.3}.\\
Notice that systems \eqref{10.1} and \eqref{10.2} are not very different from system \eqref{eq4}, the only differences are the advection terms and the homogeneity of the boundary conditions. That is why, we propose to  introduce the change of variables in \cite{meurer2012control} which allows the cancellation of the advection term.\\
To do so, we transform the boundary conditions in \eqref{10.2} to homogeneous Bcs to be able to use the transformation in \cite{meurer2012control}, we introduce the change of variables:
$$
\bar{\lambda}_{2}(x,t)=\lambda_{2}(x,t)-
(x-1)$$
which allows the homogeneity of the boundary conditions. Thus $\bar{\lambda}_{2}(x,t)$ satisfies:
\begin{equation*}
\left\{\begin{array}{llll}   
\frac{\partial }{\partial t}\bar{\lambda_{2}}(x,t)=\frac{k(t)}{C \mu} \frac{\partial}{\partial x} _{}^{C}D^{\alpha}_{x}\bar{\lambda_{2}}(x,t)-H_{1}(x)\frac{\partial }{\partial x}\bar{\lambda_{2}}(x,t)+H(x)\\
\bar{\lambda_{2}}(x,0)=\lambda_{2,0}(x)\\
\bar{\lambda_{2}}(1,t)=0\\
\frac{\partial }{\partial x}\bar{\lambda_{2}}(x,t)|_{x=0}=0.\\
\end{array}\right.
\end{equation*}
where, $H(x)=\frac{k(t)}{C \mu} \frac{\partial}{\partial x} _{}^{C}D^{\alpha}_{x}x-H_{1}(x)$ and $\lambda_{2,0}(x)=-(x-1).\\$
We introduce now the transformation in \cite{meurer2012control}:
\begin{equation} \label{meurer}  
\bar{x}=\frac{1}{\sqrt{\frac{k(t)}{c \mu}}}x
\end{equation}
the transformation in \eqref{meurer} is well defined thanks to the fact that $k(t), c$ and $\mu$ are all positive.
Thus, using transformation  \eqref{meurer},  we have:
\begin{equation} \label{meurer1}  
\lambda_{1}(x,t)=\lambda_{1}(\bar{x},t)e^{-\bar{g}_{1}(\bar{x})}
\end{equation}
and 
\begin{equation} \label{meurer2}  
\bar{\lambda}_{2}(x,t)=\bar{\lambda}_{2}(\bar{x},t)e^{-\bar{g}_{2}(\bar{x})}
\end{equation}
where, 
\begin{equation} \label{meurerr22}  
{\bar{g}_{1}(\bar{x})}=\frac{1}{2}\int_{0}^{\bar{x}} ( \frac{k(t)}{C \mu} \frac{\partial}{\partial x} _{}^{C}D^{\alpha}_{x}\lambda_{1}(x,t)-H_{1}(x)\frac{\partial }{\partial x}\lambda_{1}(x,t))|_{x=\bar{x}^{-1}(s)}ds
\end{equation}
and
\begin{equation} \label{meurerr222}  
{\bar{g}_{2}(\bar{x})}=\frac{1}{2}\int_{0}^{\bar{x}} ( \frac{k(t)}{C \mu} \frac{\partial}{\partial x} _{}^{C}D^{\alpha}_{x}\bar{\lambda_{2}}(x,t)-H_{1}(x)\frac{\partial }{\partial x}\bar{\lambda_{2}}(x,t))|_{x=\bar{x}^{-1}(s)}ds
\end{equation}
Then, systems \eqref{10.1} and \eqref{10.2} become:
 \begin{equation}\label{10.199}
\left\{\begin{array}{llll}   
\frac{\partial }{\partial t}\lambda_{1}(\bar{x},t)= \frac{\partial}{\partial \bar{x}} _{}^{C}D^{\alpha}_{\bar{x}}\lambda_{1}(\bar{x},t)+\frac{1}{C \rho}f(x)\\
\lambda_{1}(\bar{x},0)=0\\
\frac{\partial }{\partial \bar{x}}\lambda_{1}(\bar{x},t)|_{x=0}=\lambda_{1}(1,t)=0\\
\end{array}\right.
\end{equation}
and
\begin{equation}\label{10.299}
\left\{\begin{array}{llll}   
\frac{\partial }{\partial t}\bar{\lambda}_{2}(\bar{x},t)= \frac{\partial}{\partial x} _{}^{C}D^{\alpha}_{\bar{x}}\bar{\lambda}_{2}(\bar{x},t)+H(\bar{x})\\
\lambda_{2}(\bar{x},0)=\lambda_{2,0}(\bar{x})\\
\lambda_{2}(1,t)=\frac{\partial }{\partial \bar{x}}\lambda_{2}(\bar{x},t)|_{\bar{x}=0}=0.\\
\end{array}\right.
\end{equation}
Then, systems \eqref{10.199} and \eqref{10.299} can be solved following the same reasoning as the proof of Theorem \ref{thm1}. Thus:\\
\begin{equation} \label{10.16}  
\lambda_{1}(\bar{x}, t)=\frac{1}{c \rho}\int _{- \infty}^{+ \infty}  \int _{0}^{t} \mathcal{G}_{\alpha}( \bar{x}-y ,t-\tau)f(y)) d\tau dy \\
\end{equation}
where,  $\mathcal{G}_{\alpha}(\bar{x},t)$ is the Green function defined  by:
\begin{equation}\label{ggg}
\mathcal{G}_{\alpha}(\bar{x},t)=\frac{1}{2 \pi }\int_{-\infty}^{+\infty} e^{-ik\bar{x}}e^{ \phi^{\alpha+1}(s)t} dk
\end{equation}
Similarly, the solution of \eqref{10.299}:
\begin{equation} \label{10.21}  
\lambda_{2}(\bar{x}, t)=\int _{- \infty}^{+ \infty}  \mathcal{G}_{\alpha}( \bar{x}-y 
,t)\lambda_{2,0}(y) dy+(\bar{x}-1)+ \int _{- \infty}^{+ \infty}  \int _{0}^{t}   \mathcal{G}_{\alpha}( \bar{x}-y 
,t-\tau)H(y)  d\tau dy. \\
\end{equation}
The uniqueness of  \eqref{10.16} and \eqref{10.21} is a direct result from the linearity of the integral and the stability comes from  Theorem \ref{thm1} as $f(x) \in L^{1}(\mathbb{R})$. Thus :
\begin{equation} \label{m1}
\lim_{t \rightarrow +\infty} \norm{\lambda_{1}(\bar{x}, t)}_{\infty}<\infty 
\end{equation}
and 
\begin{equation} \label{m2}
\lim_{t \rightarrow +\infty} \norm{\lambda_{2}(\bar{x}, t)}_{\infty}<\infty 
\end{equation}
which is equivalent to 
\begin{equation*} 
\lim_{t \rightarrow +\infty} \norm{\lambda_{1}(x, t)}_{\infty}<\infty 
\end{equation*}
and 
\begin{equation*} 
\lim_{t \rightarrow +\infty} \norm{\lambda_{2}(x, t)}_{\infty}<\infty 
\end{equation*}
by taking the norm of \eqref{meurer1} and \eqref{meurer2} and using Cauchy–Schwarz inequality and using \eqref{m1}, \eqref{m2} and the fact that: $H_{1}(x)=\mathcal{V}^{-1}\{_{}^{C}D^{\alpha}_{y,x} R(x,y)|_{y=0}\} \in L^{1}(\mathbb{R})$.}
\bibliographystyle{plain}   
\bibliography{references.bib}   

\begin{thebibliography}{10}

\bibitem{ahmed2016adaptive}
Tarek Ahmed-Ali, Fouad Giri, Miroslav Krstic, Laurent Burlion, and
  Fran{\c{c}}oise Lamnabhi-Lagarrigue.
\newblock Adaptive boundary observer for parabolic pdes subject to domain and
  boundary parameter uncertainties.
\newblock {\em Automatica}, 72:115--122, 2016.

\bibitem{aldoghaither2017direct}
Abeer Aldoghaither, Taous-Meriem Laleg-Kirati, and Da-Yan Liu.
\newblock Direct and inverse source problems for a space fractional advection
  dispersion equation.
\newblock {\em Journal of Inverse and Ill-posed Problems}, 25(2):207--220,
  2017.

\bibitem{amir2018physics}
Sahar~Z Amir and Shuyu Sun.
\newblock Physics-preserving averaging scheme based on gr{\"u}nwald-letnikov
  formula for gas flow in fractured media.
\newblock {\em Journal of Petroleum Science and Engineering}, 163:616--639,
  2018.

\bibitem{basu2002quadratic}
Manjusri Basu and Debi~Prasad Acharya.
\newblock On quadratic fractional generalized solid bi-criterion transportation
  problem.
\newblock {\em Journal of Applied Mathematics and Computing}, 10(1-2):131,
  2002.

\bibitem{beck2012brief}
Margaret Beck.
\newblock A brief introduction to stability theory for linear pdes.
\newblock In {\em SIAM conference on Nonlinear Waves and Coherent Structures},
  2012.

\bibitem{benson2000application}
David~A Benson, Stephen~W Wheatcraft, and Mark~M Meerschaert.
\newblock Application of a fractional advection-dispersion equation.
\newblock {\em Water resources research}, 36(6):1403--1412, 2000.

\bibitem{benson1998fractional}
David~Andrew Benson.
\newblock {\em The fractional advection-dispersion equation: Development and
  application}.
\newblock PhD thesis, University of Nevada, Reno, 1998.

\bibitem{bymes2000output}
CI~Bymes, Istv{\'a}n~G Lauk{\'o}, David~S Gilliam, and Victor~I Shubov.
\newblock Output regulation for linear distributed parameter systems.
\newblock {\em IEEE Transactions on Automatic Control}, 45(12):2236--2252,
  2000.

\bibitem{byrnes1997output}
Christopher~I Byrnes, Francesco~Delli Priscoli, and Alberto Isidori.
\newblock Output regulation of nonlinear systems.
\newblock In {\em Output Regulation of Uncertain Nonlinear Systems}, pages
  27--56. Springer, 1997.

\bibitem{chen2008adi}
S~Chen and Fawang Liu.
\newblock Adi-euler and extrapolation methods for the two-dimensional
  fractional advection-dispersion equation.
\newblock {\em Journal of Applied Mathematics and Computing}, 26(1-2):295--311,
  2008.

\bibitem{deutscher2015backstepping}
Joachim Deutscher.
\newblock A backstepping approach to the output regulation of boundary
  controlled parabolic pdes.
\newblock {\em Automatica}, 57:56--64, 2015.

\bibitem{el2002continuation}
Ahmed~MA El-Sayed and Mohamed~AE Aly.
\newblock Continuation theorem of fractional order evolutionary integral
  equations.
\newblock {\em Korean Journal of Computational \& Applied Mathematics},
  9(2):525--533, 2002.

\bibitem{evans1998partial}
Lawrence~C Evans.
\newblock Partial differential equations.
\newblock {\em Graduate studies in mathematics}, 19(2), 1998.

\bibitem{fujishiro2014approximate}
Kenichi Fujishiro.
\newblock Approximate controllability for fractional diffusion equations by
  dirichlet boundary control.
\newblock {\em arXiv preprint arXiv:1404.0207}, 2014.

\bibitem{ghaffour2016class}
Lilia Ghaffour and Zoubir Dahmani.
\newblock On a class of fractional differential equations with arbitrary
  singularities.
\newblock {\em Konuralp Journal of Mathematics}, 8(2):244--251, 2016.

\bibitem{ghaffour2018fractional}
Lilia Ghaffour and Zoubir Dahmani.
\newblock Fractional differential equations with arbitrary singularities.
\newblock {\em Journal of Information and Optimization Sciences},
  39(7):1547--1565, 2018.

\bibitem{ghaffour_laleg_kirati}
Lilia Ghaffour and Taous~Meriem Laleg~Kirati.
\newblock Reference tracking problem for boundary controlled time fractional
  advection dispersion equation in the presence of disturbances.
\newblock {\em European Journal of Control, under revision}, 2021.

\bibitem{ghaffour2020non}
Lilia Ghaffour, Matti Noack, Johann Reger, and Taous-Meriem Laleg-Kirati.
\newblock Non-asymptotic state estimation of linear reaction diffusion equation
  using modulating functions.
\newblock IFAC, 2020.

\bibitem{hamalainen2010robust}
Timo H{\"a}m{\"a}l{\"a}inen and Seppo Pohjolainen.
\newblock Robust regulation of distributed parameter systems with
  infinite-dimensional exosystems.
\newblock {\em SIAM Journal on Control and Optimization}, 48(8):4846--4873,
  2010.

\bibitem{huang2005fundamental}
Fenghui Huang and Fawang Liu.
\newblock The fundamental solution of the space-time fractional
  advection-dispersion equation.
\newblock {\em Journal of Applied Mathematics and Computing}, 18(1-2):339--350,
  2005.

\bibitem{immonen2007internal}
Eero Immonen.
\newblock On the internal model structure for infinite-dimensional systems: Two
  common controller types and repetitive control.
\newblock {\em SIAM Journal on Control and Optimization}, 45(6):2065--2093,
  2007.

\bibitem{ioannou2012robust}
Petros~A Ioannou and Jing Sun.
\newblock {\em Robust adaptive control}.
\newblock Courier Corporation, 2012.

\bibitem{isidori1993topics}
A~Isidori, HW~Knobloch, and D~Flockerzi.
\newblock Topics in control theory, 1993.

\bibitem{isidori2012robust}
Alberto Isidori, Lorenzo Marconi, and Andrea Serrani.
\newblock {\em Robust autonomous guidance: an internal model approach}.
\newblock Springer Science \& Business Media, 2012.

\bibitem{jeffrey2003applied}
Alan Jeffrey.
\newblock {\em Applied partial differential equations: an introduction}.
\newblock Academic Press, 2003.

\bibitem{khalil2002nonlinear}
Hassan~K Khalil and Jessy~W Grizzle.
\newblock {\em Nonlinear systems}, volume~3.
\newblock Prentice hall Upper Saddle River, NJ, 2002.

\bibitem{krstic2008boundary}
Miroslav Krstic and Andrey Smyshlyaev.
\newblock {\em Boundary control of PDEs: A course on backstepping designs},
  volume~16.
\newblock Siam, 2008.

\bibitem{li2012observer}
C~Li, J~Wang, and J~Lu.
\newblock Observer-based robust stabilisation of a class of non-linear
  fractional-order uncertain systems: an linear matrix inequalitie approach.
\newblock {\em IET Control Theory \& Applications}, 6(18):2757--2764, 2012.

\bibitem{liang2004boundary}
Jinsong Liang, Yangquan Chen, and Rees Fullmer.
\newblock Boundary stabilization and disturbance rejection for time fractional
  order diffusion--wave equations.
\newblock {\em Nonlinear Dynamics}, 38(1-4):339--354, 2004.

\bibitem{mainardi1997fractional}
Francesco Mainardi.
\newblock Fractional calculus.
\newblock In {\em Fractals and fractional calculus in continuum mechanics},
  pages 291--348. Springer, 1997.

\bibitem{meerschaert2004finite}
Mark~M Meerschaert and Charles Tadjeran.
\newblock Finite difference approximations for fractional advection--dispersion
  flow equations.
\newblock {\em Journal of Computational and Applied Mathematics},
  172(1):65--77, 2004.

\bibitem{meurer2012control}
Thomas Meurer.
\newblock {\em Control of Higher--Dimensional PDEs: Flatness and Backstepping
  Designs}.
\newblock Springer Science \& Business Media, 2012.

\bibitem{natarajan2014state}
Vivek Natarajan, David~S Gilliam, and George Weiss.
\newblock The state feedback regulator problem for regular linear systems.
\newblock {\em IEEE Transactions on Automatic Control}, 59(10):2708--2723,
  2014.

\bibitem{paunonen2013output}
Lassi Paunonen and Seppo Pohjolainen.
\newblock Output regulation theory for distributed parameter systems with
  unbounded control and observation.
\newblock In {\em 52nd IEEE Conference on Decision and Control}, pages
  1083--1088. IEEE, 2013.

\bibitem{podlubny1998fractional}
Igor Podlubny.
\newblock {\em Fractional differential equations: an introduction to fractional
  derivatives, fractional differential equations, to methods of their solution
  and some of their applications}, volume 198.
\newblock Elsevier, 1998.

\bibitem{podlubny2007adjoint}
Igor Podlubny and YangQuan Chen.
\newblock Adjoint fractional differential expressions and operators.
\newblock In {\em Proceedings of the ASME 2007 International Design Engineering
  Technical Conferences \& Computers and Information in Engineering Conference
  IDETC/CIE}, pages 4--7, 2007.

\bibitem{schumer2001eulerian}
Rina Schumer, David~A Benson, Mark~M Meerschaert, and Stephen~W Wheatcraft.
\newblock Eulerian derivation of the fractional advection--dispersion equation.
\newblock {\em Journal of contaminant hydrology}, 48(1-2):69--88, 2001.

\bibitem{schumer2009fractional}
Rina Schumer, Mark~M Meerschaert, and Boris Baeumer.
\newblock Fractional advection-dispersion equations for modeling transport at
  the earth surface.
\newblock {\em Journal of Geophysical Research: Earth Surface}, 114(F4), 2009.

\bibitem{sheikha2009effect}
Hussain Sheikha, Mehran Pooladi-Darvish, et~al.
\newblock The effect of pressure-decline rate and pressure gradient on the
  behavior of solution-gas drive in heavy oil.
\newblock {\em SPE Reservoir Evaluation \& Engineering}, 12(03):390--398, 2009.

\bibitem{smyshlyaev2005backstepping}
Andrey Smyshlyaev and Miroslav Krstic.
\newblock Backstepping observers for a class of parabolic pdes.
\newblock {\em Systems \& Control Letters}, 54(7):613--625, 2005.

\bibitem{song2013dynamic}
Xiaona Song and Zhen Wang.
\newblock Dynamic output feedback control for fractional-order systems.
\newblock {\em Asian Journal of Control}, 15(3):834--848, 2013.

\bibitem{8483491}
Yan-Qiao Wei, Da-Yan Liu, Driss Boutat, and Hao-Ran Liu.
\newblock Robust estimation of the fractional integral and derivative of the
  pseudo-state for a class of fractional order linear systems.
\newblock In {\em 2018 37th Chinese Control Conference (CCC)}, pages
  10207--10212, 2018.

\bibitem{8866466}
Yan-Qiao Wei, Da-Yan Liu, Driss Boutat, and Hao-Ran Liu.
\newblock Non-asymptotic fractional pseudo-state differentiator for a class of
  fractional order linear systems.
\newblock In {\em 2019 Chinese Control Conference (CCC)}, pages 2347--2352,
  2019.

\bibitem{xiong2012inverse}
Xiangtuan Xiong, Qian Zhou, and YC~Hon.
\newblock An inverse problem for fractional diffusion equation in 2-dimensional
  case: Stability analysis and regularization.
\newblock {\em Journal of Mathematical Analysis and Applications},
  393(1):185--199, 2012.

\bibitem{zhou2018boundary}
Hua-Cheng Zhou and Bao-Zhu Guo.
\newblock Boundary feedback stabilization for an unstable time fractional
  reaction diffusion equation.
\newblock {\em SIAM Journal on Control and Optimization}, 56(1):75--101, 2018.

\end{thebibliography}

\end{document}